\documentclass[11pt]{amsart}

\usepackage{latexsym,amssymb}

\newtheorem{theoreme}{{\bf Th\'eor\`eme}}[section]

\newtheorem{corollaire principal}[principal]{\bf Corollaire}
\newtheorem{proposition}[theoreme]{{\bf Proposition}}
\newtheorem{lemme}[theoreme]{{\bf Lemme}}
\newtheorem{corollaire}[theoreme]{{\bf Corollaire}}

\newtheorem{remarque}{\bf Remarque}

 \newenvironment{demonstration}{\noindent{\bf D\'emonstration.
 }}{\nolinebreak $\Box $\hspace{-2.15mm}\raisebox{1.25mm}{.}  \medskip}

\newenvironment{demonstration du lemme}{\noindent{\bf D\'emonstration du lemme
}}{\nolinebreak $\Box $\hspace{-2.15mm}\raisebox{1.25mm}{.} \medskip}

\def\RR{{\bf R}}

\def\ZZ{{\bf Z}}

\begin{document}

\title[G\'eom\'etries Lorentziennes en dimension 3]{G\'eom\'etries Lorentziennes de dimension $3$ : classification et compl\'etude}

\author[S. Dumitrescu, A. Zeghib]{Sorin DUMITRESCU$^\star$ \& Abdelghani ZEGHIB$^\dagger$}

\address{${}^\star$ D\'epartement de Math\'ematique  d'Orsay, 
\'equipe de Topologie et Dynamique,
Bat. 425, U.M.R.   8628  C.N.R.S.,
Univ. Paris-Sud (11),
91405 Orsay Cedex, France}
\email{Sorin.Dumitrescu@math.u-psud.fr}
 
\address{${}^\dagger$ CNRS, UMPA, \'ecole
Normale Sup\'erieure de Lyon, France}
\email{zeghib@umpa.ens-lyon.fr}

\thanks{{\it Mots cl\'es:} Vari\'et\'es lorentziennes  localement homog\`enes- alg\`ebres de Killing transitives- compl\'etude g\'eod\'esique-compl\'etude des $(G, X)$-structures.  \\
{\it Classification math. } 53B30, 53C22, 53C50 }
 \date{le 26 octobre 2007}

\maketitle

\begin{abstract} (Lorentz geometries in dimension 3: completeness and classification). We study 3-dimensional non-Riemannian Lorentz geometries,  i.e. compact locally homogeneous Lorentz 3-manifolds, with non-compact (local)  isotropy group. One result  is that, up to a finite  cover,  all such manifolds admit  Lorentz metrics of (non-positive) constant sectional curvature. In fact, if the geometry is maximal, then, there is a tri-chotomy. The metric has constant sectionnal curvature, or is a left invariant metric on the Heisenberg group or the $SOL$-group.  These geometries, on each of the latter two groups are characterized by having  a non-compact isotropy without being flat.  Recall, for the need of his formulation of the geometrization conjecture, 
 W. Thurston  counted the 8 maximal Riemannian geometries in dimension 3.
 Here,  we count only  4 maximal  Lorentz geometries, but ignoring those  which are at the same time Riemannian. Also, all such manifolds
 are geodesically complete, except the previous non flat left invariant metric on the $SOL$-group.

 ${}$ \\

\noindent{\sc R\'esum\'e.}
Nous classifions les  g\'eom\'etries lorentziennes  de dimension $3$ qui admettent des r\'ealisations compactes. Notre r\'esultat  implique que toute vari\'et\'e lorentzienne localement homog\`ene compacte de dimension
$3$  est isom\'etrique au quotient (\`a gauche) d'un espace homog\`ene lorentzien $G/I$ de dimension $3$ par un sous-groupe discret  $\Gamma$ de $G$ agissant proprement et discontinument. Si, de plus, le groupe d'isotropie locale de la m\'etrique lorentzienne est suppos\'e non compact, le rev\^etement universel de la vari\'et\'e est isom\'etrique \`a une m\'etrique  invariante par translations sur l'un des 4 groupes suivants :
    ${\bf R}^3$, $\widetilde{SL(2, {\bf R})}$, $Heis$ ou $SOL$.
    
Notre classification implique  \'egalement que toute vari\'et\'e compacte connexe de dimension $3$ qui poss\`ede une m\'etrique lorentzienne localement homog\`ene dont le groupe d'isotropie locale est non compact, admet, \`a rev\^etement  fini pr\`es,  une m\'etrique lorentzienne de courbure sectionnelle constante (n\'egative ou nulle).

Un autre corollaire est que  toute vari\'et\'e compacte connexe de dimension $3$ munie d'une m\'etrique lorentzienne localement homog\`ene  dont le groupe d'isotropie locale est non compact est g\'eod\'esiquement compl\`ete, sauf  si elle est localement model\'ee sur l'unique m\'etrique invariante \`a gauche sur le groupe $SOL$ qui n'est pas plate, mais qui admet un groupe d'isotropie locale
non compact.

\end{abstract}



\section{Introduction, Exemples}

Soit $G$ un groupe de Lie r\'eel  et $G/I$, o\`u $I$ est un sous-groupe ferm\'e de $G$, un espace homog\`ene suppos\'e simplement connexe.

 On   dit que $G/I$
  est {\it   lorentzien} et que $(G, G/I)$ est une {\it g\'eom\'etrie lorentzienne}  (au sens de Klein) si l'action canonique de $G$ sur $G/I$ pr\'eserve une m\'etrique lorentzienne (i.e. un champ lisse de formes
  quadratiques non d\'eg\'en\'er\'ees de signature $(n-1,1)$),   ou de mani\`ere \'equivalente,  si  l'action adjointe de $I$
  pr\'eserve une forme quadratique non d\'eg\'en\'er\'ee de signature $(n-1,1)$ sur le quotient $\mathcal G / \mathcal I$ des  alg\`ebres de Lie correspondantes.

Une  vari\'et\'e   compacte connexe $M$  admet une {\it $(G, G/I)$-structure}  (on dit aussi que $M$ est {\it localement model\'ee} sur l'espace homog\`ene
  $G/I$) s'il existe un atlas de $M$ \`a valeurs dans des ouverts de $G/I$ tel que les applications de changements de cartes soient donn\'ees par des \'el\'ements du groupe $G$.
  Dans ce cas, tout objet g\'eom\'etrique (par exemple, un tenseur ou une connexion) sur $G/I$, invariant par l'action de $G$, induit un objet g\'eom\'etrique du m\^eme type sur $M$.
  En particulier, dans le cas d'un espace homog\`ene de type lorentzien  $G/I$, la vari\'et\'e $M$ h\'erite d'une m\'etrique lorentzienne localement homog\`ene (voir la d\'efinition en section~\ref{dynamique lorentzienne}). Nous dirons aussi que
  la g\'eom\'etrie lorentzienne $(G, G/I)$ admet  {\it  une r\'ealisation compacte} sur la vari\'et\'e $M$.

La   g\'eom\'etrie lorentzienne $(G, G/I)$ est dite {\it  maximale},  si l'action de $G$ ne s'\'etend pas
 en une action fid\`ele d'un groupe de Lie de dimension  strictement plus grande $G^\prime$ qui pr\'eserve une  m\'etrique lorentzienne. 
  Deux m\'etriques lorentziennes sur $G/I$ dont les  composantes  neutres des groupes des isom\'etries sont conjugu\'ees dans le groupe des diff\'eomorphismes 
  de $G/I$ d\'efinissent une m\^eme g\'eom\'etrie.

 Une $(G, G/I)$-structure sur $M$ donne classiquement naissance \`a une {\it application d\'e\-velo\-ppante} qui est un diff\'eomorphisme local entre un rev\^etement universel ${\widetilde M}$ de $M$ et
  l'espace mod\`ele $G/I$ \cite{T, Th, S}.  L'application d\'eveloppante conjugue l'action du groupe de rev\^etement ${\widetilde M} \to M$ \`a l'action du  groupe d'holonomie 
  $\Gamma$
  sur l'espace mod\`ele.

Une  $(G,G/I)$-structure est dite {\it compl\`ete} si son application d\'eveloppante est un diff\'eo\-morphisme.   
Lorsque $G$ pr\'eserve  une m\'etrique lorentzienne ou riemannienne, ou plus g\'en\'eralement une connexion, on a aussi 
la notion de {\it compl\'etude g\'eod\'esique}~\cite{Wo}.  La compl\'etude g\'eod\'esique est  plus forte  que la compl\'etude au sens ci-dessus (voir lemme~\ref{completudes}).

  Si $M$ est compl\`ete, alors $M = \Gamma \backslash G / I$, o\`u $\Gamma$ est le groupe d'holonomie de $M$, qui agit librement et   proprement sur $G/I$.  Lorsque $I$ est compact, ceci \'equivaut au fait que $\Gamma$ soit  un r\'eseau cocompact de $G$.

On dit qu'une $(G, G/I)$-structure satisfait \`a une {\it rigidit\'e de Bieberbach}, si pour toute r\'ealisation compacte compl\`ete  $M$, il existe 
un sous-groupe de Lie connexe $L$ de $G$ contenant, \`a indice fini pr\`es,  le groupe d'holonomie $\Gamma$ de $M$ et agissant librement et transitivement sur $G/I$. Dans ce cas,
$L$ s'identifie \`a $G/I$ et, \`a indice fini pr\`es, $M$ est $\Gamma \backslash L$.
 L'\'enonc\'e classique du th\'eor\`eme de Bieberbach correspond au cas de la g\'eom\'etrie euclidienne $(O(n) \ltimes {\bf R}^n, {\bf R}^n )$, avec $L = {\bf R}^n$,  le groupe des translations.

 Rappelons que  W. Thurston a classifi\'e les huit g\'eom\'etries riemanniennes  maximales  de dimension $3$ qui poss\`edent des r\'ealisations  compactes (voir \cite{T, Th, S}, et  \'egalement \cite{Dan, Souam.Toubiana},  pour des  \'etudes  r\'ecentes sur  ces g\'eom\'etries).

  Une sp\'ecificit\'e bien connue des $(G,G/I)$-structures riemanniennes est  que toute r\'ealisation compacte d'une telle structure est n\'ecessairement compl\`ete. Ce ph\'enom\`ene n'est nullement assur\'e (et souvent faux) pour les $(G,G/I)$-structures g\'en\'erales, mais sera d\'emontr\'e dans cet article pour les $(G,G/I)$-structures lorentziennes de dimension $3$.

  Nous nous int\'eressons dans ce travail aux g\'eom\'etries  lorentziennes de dimension $3$ qui 
   sont  {\it  non-riemanniennes}, i.e.    l'action de $G$ sur $G/I$ ne pr\'eserve pas de m\'etrique riemannienne. C'est \'equivalent au fait  que    l'action adjointe
  de $I$  soit  \`a image non born\'ee dans le groupe des transformations lin\'eaires de $\mathcal G  / \mathcal I$.

 Dans  la suite  de cette introduction on   explicite
des exemples de g\'eom\'etries lorentziennes non-riemanniennes maximales de dimension $3$. \\

  Les g\'eom\'etries lorentziennes de courbure sectionnelle constante sont maximales car elles r\'ealisent la dimension maximale du groupe des isom\'etries (voir proposition~\ref{dim 3,5,6} ou~\cite{Wo}). 
  
  Avant de pr\'esenter les  g\'eom\'etries mod\`eles $(G,G/I)$ qui incarnent, en dimension $3$,  les
  g\'eom\'etries lorentziennes de courbure sectionnelle constante et qui sont  la g\'eom\'etrie plate  de Minkowski (courbure nulle), la g\'eom\'etrie de Sitter (courbure positive) et la g\'eom\'etrie  anti de Sitter (courbure n\'egative), pr\'ecisons
  que, d'apr\`es un r\'esultat d\^u \`a Carri\`ere~\cite{Ca}  dans le cas plat et \'etendu  par Klingler~\cite{Kli} 
  au cas de courbure sectionnelle constante (voir \'egalement~\cite{Mess, Mor}), toute r\'ealisation compacte de $G/I$ est compl\`ete. Ce r\'esultat de compl\'etude est essentiel dans l'\'etude des g\'eom\'etries lorentziennes de courbure sectionnelle constante.

  Il implique, en particulier, que la g\'eom\'etrie de courbure sectionnelle constante positive n'admet pas de r\'ealisation compacte. En effet, d'apr\`es les travaux~\cite{CM}, 
  seuls les groupes finis agissent proprement sur l'espace de Sitter, qui n'admet donc pas de quotient compact.\\

 {\it G\'eom\'etrie Minkowski.}
  Un  mod\`ele  de la g\'eom\'etrie    Minkowski est $\RR^3$, muni de la forme quadratique $dx^2+dy^2-dz^2$. Le
  groupe $G$ de cette g\'eom\'etrie  est $O(2, 1) \ltimes \RR^3$, agissant affinement sur $\RR^3$. L'isotropie $I$ de la g\'eom\'etrie  Minkowski est $O(2, 1)$. 
  
  Il est d\'emontr\'e dans \cite{FG, GK} que la g\'eom\'etrie Minkowski satisfait  \`a une rigidit\'e de Bieberbach, avec  
   le groupe $L$  isomorphe \`a ${\bf R}^3$, $Heis$ ou $SOL$. \\

 {\it G\'eom\'etrie anti de Sitter.}
     Un mod\`ele de cette g\'eom\'etrie  est le rev\^etement universel 
 $\widetilde{SL(2, \RR)}$ de $SL(2, \RR)$, muni de la m\'etrique lorentzienne invariante par translations \`a gauche qui co\"{\i}ncide en identit\'e avec la forme de Killing $q$ sur l'alg\`ebre de Lie $sl(2, \RR)$. Comme
la forme de Killing $q$ est invariante par la repr\'esentation  adjointe, le groupe des isom\'etries de la g\'eom\'etrie  anti de Sitter contient \'egalement les translations \`a droite. La composante 
neutre  du groupe des isom\'etries de la g\'eom\'etrie  anti de Sitter est  (modulo quotient par le  noyau de l'action qui est fini, de cardinal quatre) $G= \widetilde{SL(2, \RR)} \times   \widetilde{SL(2, \RR)}$, avec isotropie $I$ isomorphe \`a $\widetilde{SL(2, \RR)}$
et   plong\'ee  diagonalement dans $G$.
  
 Comme  le groupe $\widetilde{SL(2, \RR)}$ agit librement transitivement sur le mod\`ele $G/I$, 
    il suffit de consid\'erer  le quotient \`a gauche de $\widetilde{SL(2, \RR)}$ par un r\'eseau cocompact $\Gamma$ de $\widetilde{SL(2, \RR)}$, pour construire ainsi des vari\'et\'es compactes 
    $M=\Gamma \backslash \widetilde{SL(2, \RR)}$ localement model\'ees sur la g\'eom\'etrie anti de Sitter. 
    
   La rigidit\'e de Bieberbach, valable dans le cas plat,  n'est plus valide pour  la g\'eom\'etrie  anti de Sitter~\cite{Gold,Sa}.\\

 {\it G\'eom\'etrie Lorentz-Heisenberg.} Il s'agit de la g\'eom\'etrie d'une certaine m\'etrique lorentzienne invariante par translations \`a gauche
  sur le groupe de Heisenberg $Heis$. D\'esignons par  $heis$ l'alg\`ebre de Lie de $Heis$.
    
    \begin{proposition} \label{Lorentz.Heisenberg1}  Modulo automorphisme et \`a constante multiplicative pr\`es, il existe sur $Heis$ une seule m\'etrique lorentzienne
    invariante \`a gauche, affectant une longueur positive au centre de $heis$. Ces m\'etriques d\'efinissent  une m\^eme g\'eom\'etrie lorentzienne non-riemannienne maximale, dont la composante neutre du  groupe des isom\'etries est  de dimension quatre, isomorphe \`a un produit semi-directe ${\bf R} \ltimes Heis$ et dont  l'isotropie est semi-simple  (i.e. agit sur l'espace tangent au point base comme un groupe \`a un param\`etre diagonalisable).
    Cette g\'eom\'etrie  sera  appel\'ee 
    Lorentz-Heisenberg. 
   \end{proposition}

    Comme avant, il suffit de consid\'erer un quotient (\`a gauche) de $Heis$ par un r\'eseau cocompact $\Gamma$, pour se convaincre que la g\'eom\'etrie Lorentz-Heisenberg
    se r\'ealise bien sur des vari\'et\'es compactes. \\

    {\it G\'eom\'etrie Lorentz-SOL.}  C'est une g\'eom\'etrie obtenue \`a partir d'une m\'etrique lorentzienne invariante \`a gauche sur $SOL$. Rappelons que l'alg\`ebre de Lie correspondante $sol$ est engendr\'ee  par $\{X, Z, T\}$, avec seuls crochets non-nuls $[T, X] = X$ et $[T, Z] = -Z$.  L'alg\`ebre d\'eriv\'ee est donc 
    ${\bf R}^2 = {\bf R}X \oplus {\bf R}Z$. Consid\'erons sur $SOL$ la m\'etrique lorentzienne invariante \`a gauche   $g$ d\'efinie  en identit\'e par : $X$ et $T$ sont isotropes, 
    $Z$ est orthogonal \`a ${\bf R}X \oplus {\bf R}T$, 
    et
    $g( X, T) = g( Z, Z )= 1$.

    \begin{proposition}   \label{Lorentz.SOL1}
    
 \`A automorphisme pr\`es, $g$  est l'unique m\'etrique lorentzienne invariante \`a gauche  sur $SOL$,   qui rend l'alg\`ebre d\'eriv\'ee d\'eg\'en\'er\'ee et telle que l'une des deux directions propres de $ad(T)$ est isotrope. La composante neutre du  groupe des isom\'etries de
   $g$  est de dimension $4$ : il s'agit d'une extension non triviale de $Heis$ (ici, $Heis$ n'agit pas transitivement)    et l'isotropie locale est non compacte (unipotente). 
   Cette m\'etrique d\'efinit une g\'eom\'etrie lorentzienne maximale et non-riemannienne qu'on d\'esignera par Lorentz-SOL. 
    \end{proposition}

  \section{\'Enonc\'es des r\'esultats} 
     
    Le th\'eor\`eme principal de  cet article montre que les exemples pr\'ec\'edents constituent  les seules  g\'eom\'etries lorentziennes non-riemanniennes maximales qui se r\'ealisent 
    sur des vari\'et\'es compactes de dimension $3$ :

   \begin{theoreme} \label{plein}
    Soit $M$ une vari\'et\'e lorentzienne compacte connexe de dimension $3$ localement model\'ee 
    sur une g\'eom\'etrie lorentzienne {\em non-riemannienne} $(G, G/I)$. \\
    
    Classification :  \\

    (i) Le mod\`ele $G/I$ est isom\'etrique  \`a une m\'etrique lorentzienne invariante \`a gauche sur l'un des quatre  groupes suivants~: ${\bf R}^3$, $\widetilde{SL(2, {\bf R}) }$, $Heis$ ou $SOL$.
    
       $\bullet$ Dans le cas de ${\bf R}^3$, toutes les m\'etriques sont plates.
     
     $\bullet$ Dans le cas de $\widetilde{SL(2, {\bf R}) }$,  il y a  trois g\'eom\'etries     lorentziennes non-riemanniennes qui proviennent
   des m\'etriques lorentziennes  invariantes  \`a gauche : la seule  g\'eom\'etrie maximale est anti de Sitter, les deux autres sont donn\'ees par des m\'etriques lorentziennes de courbure sectionnelle non constante  invariantes
    \`a gauche et \'egalement par un sous-groupe \`a un param\`etre unipotent  (respectivement semi-simple) de translations \`a droite.

      $\bullet$ Dans le cas de $Heis$,   il y a deux g\'eom\'etries  lorentziennes non-riemanniennes  invariantes  \`a gauche. L'une est plate  et l'autre correspond \`a la g\'eom\'etrie Lorentz-Heisenberg. 
   
  $\bullet$  Dans le  cas de $SOL$,  il y a deux g\'eom\'etries  lorentziennes non-riemanniennes invariantes \`a gauche. L'une est plate  et l'autre correspond \`a la g\'eom\'etrie Lorentz-SOL. \\
    
    (ii) Si la g\'eom\'etrie lorentzienne $(G, G/I)$ est {\em maximale,} alors  c'est l'une des 4 g\'eom\'etries suivantes : Minkowski, anti de Sitter, Lorentz-Heisenberg  ou Lorentz-SOL. \\

  Compl\'etude :\\

    (iii)  La $(G,G/I)$-structure est compl\`ete. \\ 
    
    (iv) On a une rigidit\'e de  Bieberbach dans tous les cas o\`u la g\'eom\'etrie maximale correpondante n'est pas   anti de Sitter. Dans le cas des g\'eom\'etries Lorentz-Heisenberg ou Lorentz-SOL,  le groupe d'holonomie est, \`a indice fini  pr\`es, un r\'eseau cocompact de $Heis$ ou $SOL$, respectivement.\\

    (v) $M$ est  g\'eod\'esiquement compl\`ete, sauf si elle est localement model\'ee sur  la g\'eom\'etrie Lorentz-SOL. \\
    \end{theoreme}

    Pr\'ecisons \`a pr\'esent les implications qui existent entre les diff\'erents points du th\'eor\`eme pr\'ec\'edent. Un r\'esultat  essentiel qui permet de passer 
    de la compl\'etude g\'eod\'esique \`a  la compl\'etude au sens des $(G,G/I)$-structures est le  lemme bien connu suivant~:
    
    \begin{lemme}   \label{completudes}
 
 Soit $M$ une vari\'et\'e qui admet une $(G,G/I)$-structure telle que l'action de $G$ sur $G/I$ pr\'eserve une connexion lin\'eaire $\nabla$. Si  la connexion
  induite sur $M$ par $\nabla$ est g\'eod\'esiquement compl\`ete, alors la $(G,G/I)$-structure de $M$ est compl\`ete.
 \end{lemme}    
  
  Ce r\'esultat,  combin\'e avec le th\'eor\`eme classique de compl\'etude g\'eod\'esique de Hopf-Rinow~\cite{Wo},  implique  que les $(G,G/I)$-structures riemanniennes sur des vari\'et\'es compactes sont automatiquement compl\`etes. En particulier, toute $(G,G)$-structure (o\`u le groupe de Lie $G$ agit sur lui-m\^eme par translations) sur une vari\'et\'e compacte $M$
  est compl\`ete et $M$ s'identifie au quotient de $G$ par un r\'eseau cocompact.

Gr\^ace aux r\'esultats de Carri\`ere et Klingler \cite{Car, Kli}  qui affirment que les vari\'et\'es lorentziennes compactes de courbure sectionnelle constante sont g\'eod\'esiquement compl\`etes, le lemme~\ref{completudes} fournit  \'egalement la compl\'etude des vari\'et\'es compactes localement model\'ees sur les $(G,G/I)$-structures pour lesquelles l'action de $G$ pr\'eserve une m\'etrique lorentzienne sur $G/I$ de courbure sectionnelle constante. Le th\'eor\`eme de classification des g\'eom\'etries lorentziennes non-riemanniennes maximales (partie $(ii)$ du th\'eor\`eme~\ref{plein}) permet alors de restreindre l'\'etude de la compl\'etude aux  vari\'et\'es compactes localement model\'ees sur les g\'eom\'etries Lorentz-Heisenberg
et Lorentz-SOL. La compl\'etude et la rigidit\'e de Bieberbach des r\'ealisations compactes des g\'eom\'etries Lorentz-Heisenberg et Lorentz-SOL se d\'emontre via les propositions~\ref{semi-simple} et~\ref{unipotent} respectivement. Finalement, on obtient le r\'esultat  de compl\'etude de la partie $(iii)$ 
du th\'eor\`eme principal, qui  peut s'\'enoncer \'egalement de la mani\`ere suivante :

\begin{corollaire}   \label{homogene}

Toute vari\'et\'e lorentzienne localement homog\`ene compacte connexe et de dimension $3$ est isom\'etrique 
   au quotient (\`a gauche) d'un espace homog\`ene lorentzien  $G/I$ par un sous-groupe discret $\Gamma$ de $G$ agissant proprement.
 \end{corollaire}
 
 Mentionnons  que,  dans le contexte des m\'etriques lorentziennes analytiques de dimension $3$, l'homog\'en\'eit\'e locale sur un ouvert non vide de $M$ assure l'homog\'en\'eit\'e
 locale sur $M$~\cite{Dum}  et le corollaire pr\'ec\'edent s'applique.

 Une cons\'equence du th\'eor\`eme~\ref{plein} est   le r\'esultat d'uniformisation suivant :
 
 \begin{theoreme} \footnote{Suite \`a ce travail, nous avons d\'emontr\'e  un th\'eor\`eme analogue  dans le contexte des m\'etriques riemanniennes holomorphes~ \cite{Dumizeghib}.}      \label{topologique}
 Toute vari\'et\'e compacte connexe de dimension $3$ qui poss\`ede une m\'etrique lorentzienne localement homog\`ene dont le groupe d'isotropie
locale est non compact  admet un rev\^etement fini  qui poss\`ede une m\'etrique lorentzienne de courbure sectionnelle constante (n\'egative ou nulle).
\end{theoreme}

En effet, si la g\'eom\'etrie maximale correspondante n'est pas de courbure sectionnelle constante, celle-ci co\"{\i}ncide alors avec la g\'eom\'etrie Lorentz-Heisenberg,
ou bien avec la g\'eom\'etrie Lorentz-SOL.  Dans les deux cas, la partie $(iv)$ du th\'eor\`eme~\ref{plein} assure que (\`a rev\^etement fini  pr\`es) la vari\'et\'e $M$ est un quotient (\`a gauche) de $Heis$ ou de $SOL$
par un r\'eseau. Or, aussi bien $Heis$, que $SOL$,  poss\`edent des m\'etriques lorentziennes plates invariantes \`a gauche~\cite{Rah}  qui descendent bien sur $M$.

Ainsi,  les quotients compacts de la forme $\Gamma \backslash Heis$ portent deux  g\'eom\'etries lorentziennes maximales diff\'erentes : Minkowski et Lorentz-Heisenberg.
    Ce ph\'enom\`ene est sp\'ecifique
    \`a la g\'eom\'etrie non-riemannienne : dans le contexte riemannien  il y a unicit\'e de la g\'eom\'etrie maximale port\'ee par les vari\'et\'es compactes
    de dimension trois \cite{S}.

Par ailleurs, les  parties $(i)$ et  $(iii)$ du th\'eor\`eme~\ref{plein}  permettent  de ramener l'\'etude du probl\`eme de compl\'etude g\'eod\'esique \`a la compl\'etude g\'eod\'esique des m\'etriques invariantes  \`a gauche  sur les
groupes de Lie
    ${\bf R}^3$, $\widetilde{SL(2, {\bf R}) }$, $Heis$ ou  $SOL$ et d'obtenir la partie $(vi)$ du th\'eor\`eme~\ref{plein},  en utilisant des r\'esultats connus dans ce contexte.

 Nous avons vu que la compl\'etude g\'eod\'esique est acquise dans le cas de courbure sectionnelle constante.
  La compl\'etude g\'eod\'esique  de la g\'eom\'etrie Lorentz-Heisenberg  est d\'emontr\'ee dans la section $4$ de~\cite{Rah},  o\`u le calcul explicite des g\'eod\'esiques est pr\'esent\'e. Plus g\' en\'eralement, toutes les m\' etriques lorentziennes invariantes \`a  gauche sur Heisenberg
  sont g\'eod\'esiquement compl\`etes~\cite{Gu}.
 
 Le manque de compl\'etude g\'eod\'esique de la m\'etrique Lorentz-SOL est prouv\'e dans~\cite{Gue}.
 
Finalement le cas des m\'etriques invariantes sur $\widetilde{SL(2, \RR)}$ est analys\'e dans~\cite{GL} : le r\'esultat des auteurs implique bien qu'on a compl\'etude g\'eod\'esique d\`es
que l'isotropie n'est pas compacte (voir la preuve du corollaire~\ref{sl2}). 

Par ailleurs, les exemples de~\cite{GL} et l'exemple de~\cite{Gue} montrent que la compl\'etude g\'eod\'esique peut tomber en d\'efaut pour certaines m\'etriques
lorentziennes invariantes \`a gauche sur $\widetilde{SL(2, \RR)}$ ou sur $SOL$ (ici, m\^eme en pr\'esence d'un groupe d'isotropie non compact).
N\'eanmoins notre r\'esultat de compl\'etude pour les r\'ealisations compactes des $(G, G/I)$-g\'eom\'etries  lorentziennes persiste m\^eme quand le mod\`ele $G/I$ lui-m\^eme n'est pas g\'eod\'esiquement complet.

\subsubsection*{Organisation de l'article} 
 Dans la section~\ref{dynamique lorentzienne},  nous \'etudions l'alg\`ebre de Lie des champs de Killing locaux $\mathcal G$ d'une vari\'et\'e lorentzienne localement
homog\`ene de dimension $3$.
Comme $\mathcal G$ est suppos\'ee transitive, la dimension de $\mathcal G$ est au moins $3$. Nous montrons que le cas d\'elicat auquel on peut se ramener rapidement est celui o\`u $\mathcal G$ est
de dimension $4$ et l'alg\`ebre d'isotropie locale $\mathcal I$ est de dimension $1$ (unipotente ou semi-simple). Nous construisons des feuilletages $\mathcal F$, totalement g\'eod\'esiques et d\'eg\'en\'er\'es pour
la m\'etrique, gr\^ace \`a la non-compacit\'e du groupe d'isotropie locale. 

La section~\ref{metriques invariantes} pr\'esente les g\'eom\'etries lorentziennes  invariantes par translations sur les groupes $Heis$ et $SOL$. On y d\'emontre les propositions~\ref{Lorentz.Heisenberg1} et~\ref{Lorentz.SOL1}.

Dans la  section~\ref{modeles algebriques}, on d\'etermine les g\'eom\'etries  dont l'alg\`ebre de Killing contient une copie de $sl(2, \RR)$. En particulier,
on classifie celles qui sont invariantes par translations sur $\widetilde{SL(2, \RR)}$.

Dans  la section~\ref{modeles algebriques 2},   on classifie les  alg\`ebres r\'esolubles de dimension $4$ qui peuvent appara\^{\i}tre comme
alg\`ebre de Killing d'une g\'eom\'etrie lorentzienne non-riemannienne. Une \'etape interm\'ediaire importante est la d\'etermination de la structure alg\'ebrique 
du stabilisateur d'une feuille de $\mathcal F$. On classifie ensuite les mod\`eles alg\'ebriques correspondantes $G/I$ sans hypoth\`ese  d'existence de  r\'ealisations compactes.

Lors de  la section~\ref{isotropie unipotente}, on d\'emontre que la seule g\'eom\'etrie lorentzienne non-riemanniennes maximale avec $\mathcal G$ r\'esoluble  de dimension $4$ et isotropie contenant un sous-groupe \`a un param\`etre unipotent est Lorentz-SOL. On montre \'egalement  la compl\'etude et la rigidit\'e de Bieberbach   des r\'ealisations compactes de la g\'eom\'etrie Lorentz-SOL.

Dans  la section~\ref{cas semi-simple}, on prouve que dans le cas restant, $\mathcal G$ r\'esoluble de dimension $4$ et isotropie contenant un sous-groupe \`a un param\`etre semi-simple,
la g\'eom\'etrie est Lorentz-Heisenberg. On d\'emontre ici la compl\'etude et la rigidit\'e de Bieberbach  des r\'ealisations compactes de la g\'eom\'etrie Lorentz-Heisenberg.


 \section{Dynamique lorentzienne locale}   \label{dynamique lorentzienne}

 Consid\'erons une m\'etrique lorentzienne $g$ sur une vari\'et\'e $M$. Dans toute la suite de l'article la vari\'et\'e $M$ sera suppos\'ee {\it compacte}, connexe et de dimension $3$ et $g$ sera suppos\'ee {\it localement homog\`ene}.
 
 Rappelons qu'un champ de vecteurs (local) est dit {\it champ de Killing} (local) de $g$ si son flot (local) pr\'eserve $g$. La  m\'etrique lorentzienne $g$ est dite {\it  localement homog\`ene} si l'alg\`ebre des champs de Killing locaux $\mathcal G$ de $g$ agit transitivement sur $M$. Dans ce cas, la dimension de $\mathcal G$ est minor\'ee par la dimension de $M$. La sous-alg\`ebre { \it d'isotropie} $\mathcal I$, form\'ee par les \'el\'ements de $\mathcal G$ qui s'annulent en un point $x_{0} \in M$,  engendre le pseudo-groupe des
 isom\'etries locales qui fixent $x_{0}$.
 
 Rappelons qu'une isom\'etrie locale de $g$ est enti\`erement d\'etermin\'ee par son $1$-jet~\cite{Wo, DG}, ce qui explique que $\mathcal I$ s'injecte dans le groupe
 orthogonal de $(T_{x_{0}}M, g_{x_{0}})$. Il vient que la dimension de $\mathcal I$ est inf\'erieure ou \'egale \`a la dimension du groupe orthogonal pr\'ec\'edent.

 \begin{proposition}   \label{dim 3,5,6}
 \`A rev\^etement fini pr\`es, toute m\'etrique lorentzienne localement homog\`ene $g$ sur une vari\'et\'e  compacte  $M$ de dimension $3$ est localement model\'ee  sur une unique g\'eom\'etrie maximale  $(G,G/I)$, avec $G$ connexe.
 
 i) La dimension de $G$ est $\leq 6$, avec \'egalit\'e si et seulement si $g$ est de courbure sectionnelle constante.
 
 ii) Si la dimension de $G$ est \'egale \`a $3$, alors $M$ est (\`a rev\^etement fini pr\`es) un quotient $\Gamma \backslash G$, de $G$,  par un r\'eseau cocompact $\Gamma$. De plus, l'image r\'eciproque de $g$ sur $G$ est une m\'etrique lorentzienne invariante par translation \`a gauche.
 
 iii) La dimension de $G$ est diff\'erente de $5$ (donc $I$ n'est jamais de dimension $2$).
 \end{proposition}
 
 \begin{remarque} En dimension plus grande,   en g\'en\'eral, l'existence d'un mod\`ele $G/I$ tombe en d\'efaut, aussi bien dans le contexte riemannien~\cite{K,LT}, que dans le contexte 
  pseudo-riemannien~\cite{Pa}.
  \end{remarque}
 
 \begin{demonstration}
  Un mod\`ele local $G/I$ existe  si et seulement si le groupe $I$ associ\'e \`a la sous-alg\`ebre $\mathcal I$ de $\mathcal G$ est ferm\'e dans l'unique  groupe connexe et simplement connexe $G$ associ\'e \`a $\mathcal G$ (voir le th\'eor\`eme 1.3 de~\cite{Pa}). Or, d'apr\`es un r\'esultat de G. Mostow~\cite{Mos},  ceci est vrai d\`es que la codimension
  de $\mathcal I$ dans $\mathcal G$ est $< 5$.  La codimension de $\mathcal I$ dans $\mathcal G$ \'etant ici  \'egale \`a $3$ (la dimension de $M$), le th\'eor\`eme de Mostow s'applique. 
 La vari\'et\'e $M$ est alors localement model\'ee sur une g\'eom\'etrie dont le groupe est le groupe des isom\'etries de $G/I$ (et pas seulement sa composante
  neutre $G$). Comme le  groupe des isom\'etries de $G/I$ admet un nombre fini de composantes connexes, un rev\^etement fini de $M$ est localement model\'e sur $(G,G/I).$
 
i) La dimension de $\mathcal G$ est born\'ee sup\'erieurement par la somme de la dimension de $M$ et de la dimension du groupe orthogonal correspondant. En dimension $3$,  la  dimension maximale de $\mathcal G$ est donc  \'egale \`a $6$ et elle caract\'erise  les  
 m\'etriques de courbure sectionnelle constante.  En effet, dans ce cas $\mathcal I$ est de dimension $3$ et agit transitivement sur les $2$-plans non d\'eg\'en\'er\'es contenus dans $T_{x_{0}}M$. Ceci implique, dans un premier temps, que la courbure sectionnelle en $x_{0}$ ne d\'epend que de $x_{0}$ et l'homog\'en\'eit\'e locale permet de conclure que  la courbure sectionnelle est constante sur $M$ (pour  les notions classiques de g\'eom\'etrie lorentzienne le lecteur pourra consulter~\cite{Wo}).
 
 ii) Dans ce cas $M$ admet une $(G,G)$-structure, avec $G$ agissant par translation \`a gauche sur lui-m\^eme.  Le rev\^etement universel de $M$ s'identifie donc au mod\`ele $G$ et $M$ est isom\'etrique au quotient de $G$ par un r\'eseau cocompact.
   
 iii)   La preuve est bas\'ee essentiellement sur le fait que le stabilisateur d'une orbite d'une action lin\'eaire alg\'ebrique de $PSL(2,  \RR)$ sur un espace 
  vectoriel de dimension finie est un sous-groupe alg\'ebrique qui n'est jamais de dimension $2$. Pour la preuve de ce fait il suffit de le constater pour
  les repr\'esentations irr\'eductibles de $PSL(2,\RR)$ et de remarquer que, dans le cas g\'en\'eral,  le stabilisateur d'une orbite est l'intersection des stabilisateurs
  qui correspondent aux projections de cette orbite sur chaque repr\'esentation irr\'eductible.
  
  Consid\'erons l'action de $\mathcal I$  sur l'espace tangent $T_{x_{0}}M$ (le groupe d'isotropie locale s'identifie alors \`a un sous-groupe du  groupe orthogonal
  $O(2,1) \simeq PSL(2, \RR)$) . Cette action pr\'eserve le tenseur courbure de Ricci en $x_{0}$, not\'e $Ricci_{x_{0}}$.
  
  Consid\'erons l'action du groupe orthogonal $PSL(2, \RR)$ sur l'espace vectoriel $S^2(T^*_{x_{0}}M)$ des formes quadratiques sur $T_{x_{0}}M$. L'action pr\'eserve la m\'etrique
  $g_{x_{0}}$ et induit une action de $PSL(2, \RR)$ sur l'espace vectoriel quotient $S^2(T^*_{x_{0}}M) /  \RR g_{x_{0}}$. 
  
  Le groupe d'isotropie locale est contenu dans le stabilisateur de l'\'el\'ement induit par $Ricci_{x_{0}}$ dans l'espace vectoriel $S^2(T^*_{x_{0}}M) /  \RR  g_{x_{0}}$.

  Si la dimension du groupe d'isotropie est strictement sup\'erieure \`a $1$, il vient que le stabilisateur de la courbure de Ricci est de dimension $3$ (on a vu que la dimension ne peut
  \^etre \'egale \`a $2$) et que ce stabilisateur co\"{\i}ncide avec tout le groupe orthogonal. Ceci implique que
  $Ricci_{x_{0}} = \lambda g_{x_{0}}$ avec $\lambda \in \RR$ (la fonction $\lambda$ est constante sur $M$ par homog\'en\'eit\'e locale) et notre espace est  \`a courbure sectionnelle constante. Le groupe d'isotropie est alors de dimension $3$. On vient de prouver que le groupe d'isotropie n'est jamais de dimension
  $2$ et que donc la dimension de $G$ est $\neq 5$.
 \end{demonstration}

  {\it  Il reste \`a r\'egler le cas o\`u l'alg\`ebre de Lie $\mathcal G$ des champs de Killing est de dimension $4$ et l'alg\`ebre d'isotropie $\mathcal I$ est de dimension $1$ et agit non
   proprement sur $\mathcal G$/$\mathcal I$. Ce sera toute la suite de notre travail qui se concentrera sur cette situation.}
   
     La composante neutre  du groupe d'isotropie locale  $I$ en $x_{0}$ s'identifie avec un sous-groupe \`a un param\`etre du groupe orthogonal $SO(2,1)$.
   
   Rappelons que  l'action de $SO(2,1)$ sur $(T_{x_{0}}M, g_{x_{0}})$ est conjugu\'ee   \`a l'action par la repr\'esentation  adjointe de $PSL(2, \RR)$ sur
   son alg\`ebre de Lie (car cette action pr\'eserve la forme de Killing
   sur $sl(2, \RR)$).
   
   Rappelons \'egalement que les sous-groupes \`a un param\`etre de $PSL(2, \RR)$   sont conjugu\'es \`a l'un des sous-groupes suivants :
   
   \begin{enumerate}

\item 

un sous-groupe {\it elliptique} de la forme $\left(  \begin{array}{cc}
                                                                 cos t   &   sint \\
                                                                 -sin t     &  cos t \\
                                                                 \end{array} \right)$ qui fixe un vecteur
de norme $-1$ dans $T_{x_{0}}M$;

\item   un sous-groupe {\it unipotent}  $\left(  \begin{array}{cc}
                                                                 1   &   t \\
                                                                 0     &  1\\
                                                                 \end{array} \right)$  qui fixe  un vecteur isotrope dans $T_{x_{0}}M$;

\item  un   sous-groupe {\it semi-simple}  $\left(  \begin{array}{cc}
                                                                 t   &   0\\
                                                                 0     &  t^{-1} \\
                                                                 \end{array} \right)$ qui fixe un vecteur de norme $1$ dans $T_{x_{0}}M$.

\end{enumerate}  
   
   Le cas elliptique sera exclu  de notre \'etude car il repr\'esente une g\'eom\'etrie riemannienne. Nous allons nous concentrer donc sur les cas d'isotropie semi-simple ou d'isotropie unipotente.
   
      {\bf Base adapt\'ee.} Dans le cas d'isotropie semi-simple l'action de $\mathcal I$  sur $T_{x_{0}}M$ fixe un vecteur $e_{1}$ de norme   \'egale \`a $1$. Le plan $e_{1}^{\bot}$ est alors lorentzien
   et les vecteurs $e_{2}, e_{3} \in e_{1}^{\bot}$ sont d\'efinis \`a constante multiplicative pr\`es  par les conditions suivantes : $e_{2}, e_{3}$ engendrent les deux directions isotropes du plan $e_{1}^{\bot}$
   et $g(e_{2}, e_{3})=1$. L'action du temps $t$ du flot de $\mathcal I$ s'exprime dans la base $(e_{1},e_{2}, e_{3})$ par $(e_{1}, e_{2}, e_{3}) \to (e_{1}, e^t e_{2}, e^{-t} e_{3})$.
   
   Dans le cas d'isotropie unipotente  l'action de $\mathcal I$ sur $T_{x_{0}}M$ fixe un vecteur isotrope $e_{1}$ et donc \'egalement le plan d\'eg\'en\'er\'e $e_{1}^{\bot}$ (qui contient
   $e_{1}$). 
   Consid\'erons alors des vecteurs $e_{2}$ et $e_{3}$ qui v\'erifient les relations suivantes :
    $g(e_{1}, e_{2})=0$ , $g(e_{2},e_{2})=1$, $g(e_{3},e_{3})=0,
 g(e_{2},e_{3})=0$ et $g(e_{3}, e_{1})=1.$ 
 
 Une telle base sera dite {\it adapt\'ee}. Remarquons qu'une base adapt\'ee est enti\`erement d\'etermin\'ee par le choix du vecteur unitaire $e_{2} \in e_{1}^{\bot}$~: $e_{3}$ sera
 alors l'unique vecteur isotrope situ\'e dans le plan lorentzien $e_{2}^{\bot}$ et qui est tel que $g(e_{3}, e_{1})=1$ (les vecteurs $e_{1}$ et $e_{3}$ engendrent donc les deux directions isotropes de $e_{2}^{\bot}$). Le passage d'une base adapt\'ee \`a une autre se fait par l'action de
 la diff\'erentielle en $x_{0}$ du groupe d'isotropie engendr\'e par $\mathcal I$. La matrice de cette diff\'erentielle dans la base $(e_{1}, e_{2}, e_{3})$ sera alors 
 $\left(  \begin{array}{ccc}
                                                                 1   &   t & - \frac{t^2}{2}\\
                                                                 0     &  1 &  -t\\
                                                                 0     &   0  &  1\\
                                                                 \end{array}  \right).$

   \begin{lemme}  \label{X}
   Si $\mathcal G$ est de dimension $4$, alors (\`a rev\^etement fini pr\`es) $M$ poss\`ede un champ de vecteurs  $X$ qui est $\mathcal G$-invariant et de divergence nulle.
  
  i) Si $I$ est unipotent, alors $X$ est partout    isotrope.
  
  ii) Si $I$ est semi-simple, alors $X$ est un champ de Killing de norme constante \'egale \`a $1$. 
   \end{lemme}
   
   \begin{corollaire} Si $I$ est semi-simple, $\mathcal G$ admet un centre non trivial.
   \end{corollaire}
   
    \begin{remarque} Dans le cas d'isotropie unipotente le champ de vecteurs fix\'e par l'isotropie n'est pas toujours de Killing (voir le cas $c \neq 0$ dans la section~\ref{isotropie unipotente} qui m\`ene \`a la g\'eom\'etrie Lorentz-SOL). Ce ph\'enom\`ene repr\'esente une exception lorentzienne car dans le cas des g\'eom\'etries riemanniennes
avec un groupe d'isom\'etries de dimension $4$, le champ de vecteurs stabilis\'e par l'isotropie est n\'ecessairement de Killing~\cite{Th}.
   \end{remarque}
   
  \begin{demonstration} 
   Il suffit de remarquer que le fibr\'e des rep\`eres orthonorm\'ees de $M$
   admet une r\'eduction au groupe structural $I$ (o\`u l'on identifie $I$ \`a son image dans $O(2,1)$ par sa repr\'esentation d'isotropie). Le groupe $I$ n'est pas n\'ecessairement connexe dans $O(2,1)$, mais il a un nombre fini de composantes connexes et,  \`a  rev\^etement fini pr\`es de $M$, on peut consid\'erer que $I$ est connexe
   et conjugu\'e au stabilisateur d'un vecteur de norme constante \'egale \`a $1$ (cas semi-simple) ou bien de norme constante \'egale \`a $0$ (cas unipotent). Cette  r\'eduction 
   du groupe structural d\'etermine un champ de vecteurs $X$ sur $M$, naturellement invariant par l'action de $\mathcal G$, et qui est de norme \'egale \`a $1$, si $I$ est semi-simple, ou bien isotrope, si $I$ est unipotent.

    Montrons que $X$ est de divergence nulle.  D\'esignons par 
$div (X)$ la divergence du champ de vecteurs $X$ par rapport \`a la forme volume
 $vol$ induite sur $M$ par la m\'etrique lorentzienne~: $L_{X} vol =div(X) \cdot vol,$ o\`u $L_{X}$ est la d\'eriv\'ee de Lie dans la direction de $X$. Comme $X$ et $vol$ sont
 $\mathcal G$-invariants, il vient que la fonction $div(X)$ est \'egalement $\mathcal G$-invariante et donc constante \'egale \`a $\lambda \in \RR$.
Notons par $\phi^t$ le temps $t$ du flot de $X$ : nous avons alors que $(\phi^t)^*vol=exp( \lambda  t) \cdot vol,$ pour tout $t \in \RR$. Comme
le flot de $X$ doit pr\'eserver $\int_{M}vol$, il vient que $\lambda =0$.  

ii) Supposons que $I$ est semi-simple et montrons que $X$ est de Killing.   On montre d'abord qu'\`a l'instar de  $\mathcal G$, l'action de $\phi^t$  pr\'eserve  $X^{\bot}$. En effet, fixons un point $x_{0} \in M$ et consid\'erons son image $\phi^t(x_{0})$ par le temps $t$ du flot de $X$. Pour chaque $t$ consid\'erons
 une isom\'etrie locale $g^t$ qui envoie $x_{0}$ sur $\phi^t(x_{0})$.
    
    Les diff\'eomorphismes locaux  $(g^t)^{-1} \circ \phi^t$ fixent $x_{0}$ ainsi que le vecteur $X(x_{0}) \in T_{x_{0}}M$,  et commutent avec toutes les isom\'etries locales.  Donc,  la diff\'erentielle $L_{t}$ de $(g^t)^{-1} \circ \phi^t$
    en $x_{0}$ commute avec l'action de l'isotropie en $x_{0}$ et pr\'eserve, par cons\'equent,  les sous-espaces vectoriels stabilis\'es par l'isotropie.
    
    La diff\'erentielle $L_{t}$ pr\'eserve donc le plan $X(x_{0})^{\bot}$ et ses deux  droites isotropes. Comme $div(X)=0$, la diff\'erentielle
    $L_{t}$ pr\'eserve le volume. Il  vient que   le produit des  deux valeurs propres associ\'ees aux deux  directions isotropes de $X(x_{0})^{\bot}$ vaut  $1$. Ceci implique
    que la diff\'erentielle de      $(g^t)^{-1} \circ \phi^t$ en $x_{0}$ est une isom\'etrie et que le flot de $X$ agit par isom\'etries.
    
    Nous avons donc que $X$ est un \'el\'ement non trivial du centre de $\mathcal G$.                           
     \end{demonstration}                 
                                                            
          D\'esignons par $\nabla$ la connexion de Levi-Civita de $g$. Nous avons le 
                                                                 
           \begin{lemme}   \label{invariance}
     
    Si $I$ est unipotent, alors   le champ d'endomorphismes $\nabla_{\cdot}X$ du fibr\'e tangent     s'exprime dans une base adapt\'ee sous la forme 
 $ \left(  \begin{array}{ccc}
                                                                 0    &   0 & \alpha \\
                                                                 0     &  0  &  0 \\
                                                                 0     &   0  &  0 \\
                                                                 \end{array} \right),$ avec $\alpha \in \RR$.
                                                                 
   De plus, $X$ est de Killing si et seulement si $\alpha =0.$
                                                                 
      \end{lemme}
      
      \begin{demonstration}
      Soit $B$ la matrice de  $\nabla_{\cdot}X$  dans une base adapt\'ee de $T_{x_{0}}M$ .
      L'invariance de $\nabla_{\cdot}X$   par l'action de $\mathcal I$ implique, en particulier, que la  matrice $B$ commute avec $L_{t}=\left(  \begin{array}{ccc}
                                                                 1   &   t & - \frac{t^2}{2}\\
                                                                 0     &  1 &  -t\\
                                                                 0     &   0  &  1\\
                                                                 \end{array}  \right),$  pour tout $t \in \RR$.

      Les espaces propres de chacune des deux matrices sont stables par l'autre. Comme $L_{t}$ ne laisse stable aucune d\'ecomposition
      non triviale de $T_{x_{0}}    M$ en somme directe, il vient que toutes les valeurs propres de $B$ sont \'egales.               Un calcul direct pr\'ecise la forme de la matrice~: 
      $ B=  \left(  \begin{array}{ccc}
                                                                 \lambda    &   \beta  & \alpha \\
                                                                 0     &  \lambda  &  - \beta  \\
                                                                 0     &   0  &  \lambda \\
                                                                 \end{array} \right),     \hbox{avec} \; \;   \alpha, \beta, \gamma \in \RR.$

 Il suffit de montrer que $\nabla_{\cdot }     X$  est nilpotent et que $\beta =0$. 
 
 D'apr\`es le lemme~\ref{X}, $div(X)=0$. 
Or, la   divergence de $X$ en un point $x_{0}$ de $M$ n'est rien d'autre que la trace de l'endomorphisme
$(\nabla_{\cdot} X)(x_{0})$, \'egale en occurrence \`a $3 \lambda$. 

Ainsi, $\lambda =0$. Par ailleurs, il sera montr\'e, de mani\`ere ind\'ependante,  dans  la proposition~\ref{feuilles plates},  que $X$ est parall\`ele le long de toute courbe tangente \`a $X^{\bot}$. Ceci implique que $\nabla_{e_{2}} X=0$ et donc que $\beta =0$.

Le champ de vecteurs $X$ est de Killing si et seulement si l'op\'erateur $\nabla_{\cdot} X$ est $g$-anti-sym\'etrique~\cite{Wo}. Un op\'erateur de rang $\leq 1$ \'etant  anti-sym\'etrique
si et seulement s'il est identiquement nul, il vient que $X$ est de Killing si et seulement si $\alpha =0.$
    \end{demonstration}

        Un ph\'enom\`ene caract\'eristique \`a la dynamique lorentzienne, remarqu\'e pour la premi\`ere fois par M. Gromov dans~\cite{Gro} (voir \'egalement l'article de survol~\cite{DG})
        et qui a \'et\'e amplement utilis\'e depuis est le fait que la pr\'esence d'un groupe d'isotropie (local) non compact implique l'existence de hypersurfaces
        totalement g\'eod\'esiques d\'eg\'en\'er\'ee pour la m\'etrique. 
        
       Le fait remarquable est que dans notre contexte ce feuilletage est r\'egulier : 
        
        \begin{lemme}   \label{feuilletages}

        i) Dans le cas d'isotropie unipotente la vari\'et\'e $M$ poss\`ede  un feuilletage $\mathcal F$ de dimension deux,  totalement g\'eod\'esique et $g$-d\'eg\'en\'er\'e,  dont le plan tangent en chaque point
        est $X^{\bot}.$
        
        ii)  Dans le cas d'isotropie semi-simple la vari\'et\'e $M$ poss\`ede deux feuilletages de dimension deux, totalement g\'eod\'esiques et d\'eg\'en\'er\'ees $\mathcal F$$_{1}$ et $\mathcal F$$_{2}$.
        Le plan tangent \`a chacun de ces feuilletages est engendr\'e par le champ de vecteurs $X$   et par une des deux droites
        isotropes de $X^{\bot}$.
        
        L'action de $\mathcal G$ pr\'eserve chacun de ces feuilletages.
        \end{lemme}
        
        \begin{demonstration}
        L'id\'ee de la d\'emonstration consiste \`a voir le graphe d'un  \'el\'ement $f$ du groupe d'isotropie locale en $x_{0} \in M$ comme  une sous-vari\'et\'e de dimension $3$
        au voisinage de $(x_{0} , x_{0})$ dans $M \times M$. Cette sous-vari\'et\'e est totalement g\'eod\'esique et isotrope pour la m\'etrique pseudo-riemannienne $g \oplus(-g)$
        sur $M \times M$. Si $f_{n}$ est une suite d'\'el\'ements du groupe d'isotropie locale qui tend vers l'infini dans $O(2,1)$, alors la suite des graphes de $f_{n}$ converge vers une sous-vari\'et\'e
        $F'$ totalement g\'eod\'esique et isotrope, mais qui ne repr\'esente plus un graphe. N\'eanmoins l'intersection de $F'$ avec l'espace vertical $\{x_{0} \} \times M$ est
        de dimension au plus $1$,  car il s'agit d'une sous-vari\'et\'e isotrope  de la vari\'et\'e lorentzienne $M$. Il vient que la projection $F$ de $F'$ sur l'horizontale
        $M \times \{x_{0} \}$ s'identifie avec une surface totalement g\'eod\'esique qui passe par $x_{0}$ et qui est d\'eg\'en\'er\'ee. Rappelons ici que le rang et la signature de la
        restriction de $g$ \`a une surface totalement g\'eod\'esique sont invariants par transport parall\`ele et donc constants.
        
        Dans notre cas,  il suffit de consid\'erer  une suite d'isom\'etries locales qui se trouve dans le sous-groupe \`a un param\`etre engendr\'e par $\mathcal I$. Ces isom\'etries
        se lin\'earisent en coordonn\'ees exponentielles et s'expriment dans une base adapt\'ee de $T_{x_{0}}M$ sous la forme pr\'esent\'ee pr\'ec\'edemment. On constate
        imm\'ediatement que la limite de nos suites de graphes d'applications lin\'eaires est  le plan $X(x_{0})^{\bot}$ dans le cas d'isotropie unipotente
        et les deux plans engendr\'es par $X(x_{0})$ et par chacune des deux directions isotropes de $X(x_{0})^{\bot}$ dans le cas d'isotropie semi-simple.
        
        Finissons la preuve dans le cas d'isotropie unipotente.
        
        Nous venons de prouver que par chaque point $x_{0} \in M$ passe une surface totalement g\'eod\'esique d\'eg\'en\'er\'ee tangente \`a $X(x_{0})^{\bot}$. Nous d\'emontrons que cette surface est unique. Supposons   par l'absurde qu'il existe une deuxi\`eme surface totalement g\'eod\'esique d\'eg\'en\'er\'ee tangente au point $x_{0}$ \`a un plan d\'eg\'en\'er\'e contenu
        dans $T_{x_{0}}M$. Comme $\mathcal I$ agit transitivement sur les plans d\'eg\'en\'er\'es diff\'erents de $X(x_{0})^{\bot}$, il vient que chaque plan d\'eg\'en\'er\'e contenu
        dans  $T_{x_{0}}M$ est tangent \`a une surface totalement g\'eod\'esique    d\'eg\'en\'er\'ee. Ceci implique que la m\'etrique lorentzienne $g$ est \`a courbure sectionnelle
        constante (voir~\cite{Z}, proposition 3)~: absurde (car nous sommes dans le cas d'isotropie de dimension $1$).
        
        Il vient que l'unique surface totalement g\'eod\'esique d\'eg\'en\'er\'ee qui passe par $x_{0}$ est en tout point tangente  \`a $X^{\bot}$. Le champ de plans $X^{\bot}$
        est int\'egrable et le feuilletage $\mathcal F$ engendr\'e est totalement g\'eod\'esique et d\'eg\'en\'er\'e.
        
        Comme l'action de $\mathcal G$ pr\'eserve $X$, elle pr\'eserve \'egalement $X^{\bot}$ et le feuilletage $\mathcal F$.
        Les m\^emes arguments s'appliquent dans le cas d'isotropie semi-simple.
         \end{demonstration}

  {\bf Le stabilisateur $H$ d'une feuille}. Si $I$ est unipotent, d\'esignons par $\mathcal H$ la sous-alg\`ebre de $\mathcal G$ qui stabilise la feuille $F=\mathcal F$$(x_{0})$ de $x_{0}$ et par $H$ le sous-groupe de $G$ correspondant. Si $I$ est semi-simple, on garde les m\^emes notations pour le stabilisateur de $\mathcal F_{1}$$(x_{0})$.
  
  \begin{lemme} Le groupe $H$ est de dimension $3$ et agit transitivement sur $\mathcal F$$(x_{0})$ (respectivement $\mathcal F_{1}$$(x_{0})$). L'isotropie $I$ en $x_{0}$ est
  contenue dans $H$.
  \end{lemme}
  
  \begin{corollaire} Les feuilles de $\mathcal F$ (resp. $\mathcal F_{1}$) sont localement model\'ees sur $(H,H/I)$.
  \end{corollaire}
  
  \begin{demonstration}
  Faisons la preuve dans le cas o\`u $I$ est unipotent.
  
 Une isom\'etrie locale de $g$ qui relie deux points $x_{0}$ et $x_{1}$ d'une m\^eme feuille $\mathcal F$$(x_{0})$ pr\'eserve $X$ et  \'egalement $X^{\bot}$ et envoie
 donc $exp_{x_{0}}(X^{\bot})$ sur $exp_{x_{1}}(X^{\bot})$. Il vient que cette isom\'etrie locale stabilise la feuille totalement g\'eod\'esique $\mathcal F$$(x_{0})$ et que, en particulier, l'alg\`ebre  d'isotropie
 $\mathcal I$ du point $x_{0}$ est contenue dans l'alg\`ebre   $\mathcal H$ qui stabilise  la feuille $\mathcal F$$(x_{0})$. Comme $\mathcal G$ agit transitivement sur $M$, la remarque pr\'ec\'edente implique
 que  $\mathcal H$  agit transitivement sur  $\mathcal F$$(x_{0})$ (avec une isotropie de dimension $1$) ce qui implique que $\mathcal H$ est de dimension $3$.
  \end{demonstration}

 Remarquons  que la restriction de $\mathcal H$ \`a une feuille $F$ de $\mathcal F$ est un isomorphisme d'alg\`ebres de Lie. Ceci est d\^u  au fait
 qu'une isom\'etrie locale de $g$ est compl\'et\'ement d\'etermin\'ee par sa restriction \`a un sous-ensemble invariant totalement g\'eod\'esique de dimension $2$ (il suffit de lin\'eariser l'isom\'etrie en coordonn\'ees exponentielles et de  v\'erifier ce fait pour les applications lin\'eaires).
 
   Dans la suite de l'article nous allons \'egalement d\'esigner par $X$ et par $\mathcal F$, les objets correspondants sur le mod\`ele $G/I$.

  \section{G\'eom\'etrie des m\'etriques invariantes \`a gauche}   \label{metriques invariantes}
   
   Dans cette section,  nous examinons la g\'eom\'etrie de certaines  m\'etriques lorentziennes invariantes par translations \`a gauche sur des groupes de Lie unimodulaires de
   dimension $3$.  Cette  \'etude a \'et\'e initi\'ee par J. Milnor~\cite{Mi}  dans le cas riemannien et poursuivie dans le cadre lorentzien \cite{No, Ra, CP}. Comme les m\'etriques
   lorentziennes qui nous int\'eressent sont pr\'ecis\'ement celles de courbure sectionnelle non constante et dont le groupe d'isotropie locale est non compact, les groupes
   $\RR^3$, $S^3$ et le groupe des d\'eplacements du plan euclidien seront exclus de notre \'etude. Il reste \`a examiner les groupes $Heis$, $SOL$ et $\widetilde{SL(2, {\bf R})}$.  
   Le cas $\widetilde{SL(2, {\bf R})}$ fera  l'objet de la proposition~\ref{semi-simple algebrique}, en  section~\ref{modeles algebriques}.

\subsection{Cas du groupe $Heis$} 

Nous d\'ecrivons ici les m\'etriques lorentziennes invariantes \`a gauche sur Heisenberg et d\'emontrons au passage la proposition~\ref{Lorentz.Heisenberg1}. \\

\begin{demonstration}
L'alg\`ebre de Lie $heis$ est engendr\'ee par un \'el\'ement central $X'$ et par deux \'el\'ements $Z$ et $T$ tels que
$\lbrack Z, T \rbrack =X'$. Les automorphismes de l'alg\`ebre de Lie pr\'eservent le centre et envoient donc l'\'el\'ement $X'$ sur un multiple de la forme $\lambda X'$, avec $\lambda \in \RR^{*}$. Un tel automorphisme agit  sur le plan $heis/ \RR X' \simeq \RR^2$ par un automorphisme de d\'eterminant \'egal \`a $\lambda$. Inversement, toutes ces transformations
sont bien des automorphismes de l'alg\`ebre de Lie $heis$.

Il est montr\'e dans~\cite{Ra} et \cite{Rah}  que, modulo automorphisme de $heis$,  il existe trois  classes de m\'etriques lorentziennes invariantes sur $Heis$
selon que l'\'el\'ement $X'$ est de norme  strictement  n\'egative, nulle  ou bien strictement positive. Pour les r\'esultats suivants, qui d\'ecrivent ces  m\'etriques, nous renvoyons
le lecteur \`a~\cite{Ra} et~\cite{Rah}.

- Dans le cas o\`u $X'$ est de norme nulle,  la m\'etrique est plate.

- Quand $X'$ est de norme constante \'egale \`a $-1$, son orthogonal est le plan riemannien engendr\'e par $Z$ et $T$ (\`a automorphisme pr\`es) et le sous-groupe \`a un param\`etre d'automorphismes
de $heis$ qui stabilise $X'$ et agit par rotations euclidiennes sur le plan engendr\'e par $Z$ et $T$ constitue un sous-groupe \`a un param\`etre (elliptique) d'isom\'etries.
Le groupe des isom\'etries  est dans ce cas de dimension $4$ et sa composante connexe de l'identit\'e est engendr\'ee donc par les translations \`a gauche et par le sous-groupe
\`a un param\`etre (d'isotropie) pr\'ec\'edent. Il vient que cet exemple est riemannien.

- Le cas int\'eressant  \`a notre sens et qui fournit la g\'eom\'etrie Lorentz-Heisenberg est celui o\`u $X'$ est de norme constante positive. Dans ce cas, quitte \`a appliquer un automorphisme de $heis$ on peut supposer que  $X'^{\bot}$  est le plan lorentzien engendr\'e par les droites isotropes port\'ees 
par $Z$ et $T$. En appliquant finalement  un automorphisme  de la forme $(X',Z,T) \to (\lambda^2 X',\lambda Z, \lambda T)$, avec $\lambda \in \RR^*$,
on peut rendre $g$ multiple de l'unique m\'etrique qui attribue \`a $X'$ la  norme $1$ et qui est telle que  $g(Z,T)=1$.

Le groupe des isom\'etries de cette famille de  m\'etriques est le m\^eme. Il est de dimension $4$, il contient $Heis$  et l'isotropie contient le sous-groupe \`a un param\`etre d'automorphismes de $heis$ qui fixe
$X'$ et qui agit sur le plan lorentzien engendr\'e par $Z$ et $T$ par les matrices  $\left(  \begin{array}{cc}
                                                                 e^t   &   0\\
                                                                 0     &  e^{-t} \\
                                                                 \end{array} \right)$. Le groupe d'isotropie contient alors  un sous-groupe \`a un param\`etre semi-simple et on est en pr\'esence
                                                                 d'une g\'eom\'etrie lorentzienne non-riemannienne dont la composante neutre du  groupe des isom\'etries est $G = \RR \ltimes Heis$, o\`u l'action du  facteur $\RR$ qui engendre
   l'isotropie sur $Heis$ vient d'\^etre explicit\'ee.

   Comme la composante neutre du  groupe des isom\'etries de cette m\'etrique n'est pas contenu dans le groupe des isom\'etries de la m\'etrique  plate~\cite{Rah},  il vient que cette g\'eom\'etrie est \'egalement maximale.
 \end{demonstration}

  \subsection{Cas du groupe $SOL$}
 Rappelons que l'alg\`ebre de Lie $sol$ est engendr\'ee  par $\{X', Z, T\}$, avec les \'el\'ements  $X'$ et $Z$ qui commutent et les deux crochets non-nuls $\lbrack T, X' \rbrack = X'$ et 
  $\lbrack T, Z \rbrack= -Z$.  
  
  Nous d\'emontrons  \`a
      pr\'esent la proposition~\ref{Lorentz.SOL1}. 
  
  \begin{demonstration}  Les automorphismes de $sol$ agissent sur l'alg\`ebre d\'eriv\'ee ${\bf R}^2 = {\bf R}X' \oplus {\bf R}Z$ en pr\'eservant la d\'ecomposition en espaces propres de l'op\'erateur
    $ad(T)$. Cette action se fait donc par homoth\'etie sur chacune des droites ${\bf R}X'$ et    $ {\bf R}Y$. Par ailleurs, un automorphisme de $sol$ envoie $T$ sur
    la somme de $T$ avec un \'el\'ement de       l'alg\`ebre d\'eriv\'ee. Inversement tout isomorphisme de l'espace vectoriel $sol$ de la forme pr\'ec\'edente  est  
  un automorphisme de $sol$.
  
  Consid\'erons une m\'etrique lorentzienne $g$,  invariante \`a gauche sur $SOL$,  qui rend l'alg\`ebre d\'eriv\'ee d\'eg\'en\'er\'ee et l'une des deux directions propres
    de $ad(T)$  isotrope.
  
    Pour fixer les id\'ees, supposons que $X'$ engendre l'unique direction isotrope de l'alg\`ebre d\'eriv\'ee et donc que ${\bf R}X' \oplus {\bf R}Z$  co\"{\i}ncide avec ${X'}^{\bot}$.
    Quitte \`a appliquer une homoth\'etie sur la direction $\RR Z$,  on peut supposer que $Z$ est de norme \'egale \`a $1$. En additionnant \`a $T$ un multiple
    de $Z$, on peut consid\'erer que $T$ est orthogonal \`a $Z$. On applique ensuite une homoth\'etie sur la droite engendr\'ee par $X'$ de mani\`ere \`a avoir $g(T,X')=1$
    (ceci est possible car le plan $Z^{\bot}$, engendr\'e par $X'$ et $T$, est lorentzien et donc $g(T,X')$ est non nul). Finalement,  on ajoute \`a $T$ un multiple de $X'$, de mani\`ere
    \`a rendre $T$ de norme \'egale \`a $0$ (tout en pr\'eservant l'orthogonalit\'e entre $T$ et $Z$). Nous venons de prouver l'unicit\'e de $g$ (\`a automorphisme pr\`es).
    
    Nous montrons \`a pr\'esent que cette m\'etrique n'est pas plate. Un calcul direct implique que $\nabla_{Z}T=Z$, $\nabla_{T}Z=0$, $\nabla_{T}T=-T$, $\nabla_{X'}X'=0$,
    $\nabla_{Z}Z=-X'$, $\nabla_{X'}Z=\nabla_{Z}X'=0$, $\nabla_{X'}T=0$ et $\nabla_{T}X'=X'$.
    
    Ceci donne que $R(T,Z)$ est un endomorphisme non nul qui s'exprime dans la base $(X',Z,T)$ par la matrice $ \left(  \begin{array}{ccc}
                                                                 0  &   -2 &  0\\
                                                                 0     &  0 & 2\\
                                                                 0     &   0  &  0\\
                                                                 \end{array}  \right). $    La m\'etrique
                                                                 n'est donc pas plate.

          Pr\'ecisons, par ailleurs, que la m\'etrique \'etant isom\'etrique \`a ses multiples par des constantes, tous ses invariants scalaires sont nuls. En particulier,
          ses courbures principales sont nulles.            La structure locale des  m\'etriques lorentziennes localement homog\`enes non plates dont tous les invariants scalaires sont nuls a \'et\'e classifi\'ee par V. Patrangenaru (voir la proposition
    3.2 de~\cite{Pa}) :  parmi ces exemples, seule la m\'etrique Lorentz-SOL se r\'ealise sur des vari\'et\'es compactes.\\

        Maintenant nous prouvons que le groupe des isom\'etries de la g\'eom\'etrie Lorentz-SOL est de dimension $4$. Pour cela nous construisons la m\'etrique
        pr\'ec\'edente d'une mani\`ere diff\'erente.     Consid\'erons l'alg\`ebre de Lie $heis$ engendr\'ee par trois \'el\'ements $Y,X',Z$, avec $X'$ \'el\'ement central
        et $\lbrack Y, Z \rbrack =X'$. Consid\'erons l'action d'un quatri\`eme \'el\'ement  $T$ sur     $heis$ donn\'ee, dans la base $(X',Z,Y)$, par la matrice                                                 
 $ad(T)= \left(  \begin{array}{ccc}
       1  &   0  & 0 \\
       0     &  -1 & 1 \\
       0 & 0 & 2
  \end{array} \right).$

  Cette  action est  bien une d\'erivation, ce qui implique que les \'el\'ements $\{ T,X',Z,Y  \}$ engendrent une alg\`ebre de Lie ${\mathcal G}$ isomorphe \`a l'alg\`ebre de Lie du produit semi-direct $G= \RR \ltimes Heis$, o\`u l'action du facteur $\RR$, engendr\'e par $T$,  sur $Heis$ vient d'\^etre explicit\'ee. Le sous-groupe \`a un param\`etre $I$ engendr\'e par $Y$ est ferm\'e dans $G$ (voir le lemme~\ref{dim 3,5,6}). L'espace homog\`ene $G/I$ est
  lorentzien car l'action de $ad(Y)$ sur ${\mathcal G} / {\mathcal I}$ s'exprime, dans 
la base $(X', Z, T)$,  par la matrice 
$ad(Y)=\left(  \begin{array}{ccc}
        0 &   1  & 0 \\
       0     &  0 & -1 \\
       0 & 0 &  0
  \end{array} \right),$ qui engendre  un sous-groupe \`a un param\`etre unipotent pr\'eservant une m\'etrique lorentzienne pour laquelle  $X'$ est isotrope et le plan engendr\'e par $X'$ et $Z$ d\'eg\'en\'er\'e. Par ailleurs, l'isotropie $I$ est
  engendr\'ee par $\RR Y$ qui intersecte trivialement l'alg\`ebre de Lie $sol$ engendr\'ee par $\{ T,X',Z  \}$. Il vient que l'action du groupe $SOL$ sur $G/I$  est libre et transitive
  et la m\'etrique pr\'ec\'edente s'identifie \`a l'unique m\'etrique invariante \`a gauche sur $SOL$ qui rend $X'$ isotrope et l'alg\`ebre d\'eriv\'ee  $\RR X' \oplus \RR Z$ d\'eg\'en\'er\'ee.
  On vient de d\'emontrer que cette m\'etrique admet une isotropie non triviale unipotente et que son groupe des isom\'etries est de dimension au moins $4$. La g\'eom\'etrie Lorentz-SOL est donc non-riemannienne.
 Comme la
  m\'etrique n'est pas de courbure sectionnelle constante, le lemme~\ref{dim 3,5,6} implique que le groupe des isom\'etries est de dimension \'egale \`a $4$ et
  sa composante neutre  co\"{\i}ncide avec $G=  \RR \ltimes Heis$. Mentionnons que, contrairement \`a la g\'eom\'etrie Lorentz-Heisenberg, ici le facteur $Heis$ contient l'isotropie $\RR Y$ et n'agit donc
  pas librement.
  
  Le groupe des isom\'etries de la g\'eom\'etrie Lorentz-SOL n'est  contenu strictement dans le groupe des isom\'etries d'aucune autre m\'etrique  invariante sur $SOL$. En effet, le groupe d'isotropie $I$ agit
  de mani\`ere unipotente sur le plan engendr\'e par $X'$ et $Z$. Ceci implique que toute m\'etrique lorentzienne invariante \`a gauche sur $SOL$, qui est pr\'eserv\'ee par $I$, rend l'alg\`ebre 
  d\'eriv\'ee d\'eg\'en\'er\'ee et la direction $\RR X'$ isotrope. Or, on a vu que ces m\'etriques co\"{\i}ncident n\'ecessairement (\`a automorphisme pr\`es) avec la  m\'etrique Lorentz-SOL.
  La g\'eom\'etrie Lorentz-SOL est donc maximale.
   \end{demonstration}
  
 Le groupe $SOL$ admet \'egalement des m\'etriques plates invariantes par translations. Pour s'en convaincre,  il suffit d'exhiber une copie  de  $SOL$ dans  le groupe des isom\'etries
   de Minkowski, qui agit simplement transitivement. Nous ne construirons pas un tel exemple ici, car  une m\'etrique plate explicite sur $SOL$ appara\^{\i}tra naturellement en section~\ref{cas semi-simple}, au  cours de la preuve de la proposition~\ref{non maximale}.

   \section{Mod\`eles alg\'ebriques $SL(2, \RR)$-invariants}  \label{modeles algebriques}

   Nous classifions ici les g\'eom\'etries  lorentziennes   non-riemanniennes de la forme $G/I$, avec $G$ groupe de Lie de dimension $4$  non r\'esoluble et  isotropie $I$ de dimension $1$.

   \begin{proposition}   \label{semi-simple algebrique}
   
   Si $\mathcal G$ a une partie semi-simple non triviale, alors $\mathcal G$ est l'alg\`ebre de Lie de $\RR \times \widetilde{SL(2, \RR)}$ et $g$ est  isom\'etrique \`a :
   
   (1)   une m\'etrique lorentzienne
   sur $\widetilde{SL(2, \RR)}$ invariante par les translations \`a gauche et par un sous-groupe \`a un param\`etre unipotent (respectivement semi-simple) de translations \`a droite. L'isotropie $I$
   est le premier facteur du produit $\RR \times \widetilde{SL(2, \RR)}$. La g\'eom\'etrie maximale correspondante  est celle de courbure sectionnelle constante n\'egative donn\'ee par la forme de Killing.
   
   (2) la g\'eom\'etrie produit de  $\RR$ et   du  plan de Sitter de dimension 2, i.e. la  surface  lorentzienne homog\`ene simplement connexe  de courbure sectionnelle constante non nulle. L'isotropie $I$ est incluse dans
   le facteur $ \widetilde{SL(2, \RR)}$   et est conjugu\'ee \`a un  sous-groupe \`a un param\`etre semi-simple.
    \end{proposition}
    
    \begin{remarque} On montrera \`a la proposition~\ref{non realisation compacte} (section~\ref{cas semi-simple})  que la g\'eom\'etrie produit du cas (2) ne se r\'ealise pas sur des vari\'et\'es compactes.
    \end{remarque}
    
    \begin{corollaire}    \label{sl2}
   
   Toute r\'ealisation compacte d'une g\'eom\'etrie  du type pr\'ec\'edent admet un rev\^etement universel isom\'etrique \`a 
$\widetilde{SL(2, \RR)}$ muni d'une m\'etrique invariante \`a gauche et par un sous-groupe \`a un param\`etre (semi-simple ou unipotent)  de translations \`a droite. Ces m\'etriques sont g\'eod\'esiquement compl\`etes.
 \end{corollaire}   
   
   Passons \`a la preuve de la proposition.

   \begin{demonstration}  Comme $\mathcal G$ est de dimension $4$ et qu'il n'existe pas de groupe de Lie semi-simple de dimension $4$, il vient
  que $\mathcal G= \RR \oplus {\mathcal G}_{1}$, o\`u  ${\mathcal G}_{1}$ est une alg\`ebre de Lie semi-simple de dimension $3$~\cite{Kir}.  L'alg\`ebre de Lie $\mathcal G$$_{1}$ est alors isomorphe \`a  $sl(2, \RR)$ ou \`a l'alg\`ebre de Lie de $S^3$,
  o\`u $S^3$ est la sph\`ere de dimension $3$ avec sa structure canonique de groupe de Lie.\\

  (1) Consid\'erons d'abord  le cas o\`u $\mathcal G$$_{1}$ agit librement transitivement sur $G/I$ et donc $g$ s'identifie avec une m\'etrique invariante par translations sur le groupe de
    Lie associ\'e $G_{1}$. 
  
  Si $G_{1}$ est
  $S^3$,  le mod\`ele $G/I$ s'identifie \`a  $S^3$ muni d'une m\'etrique lorentzienne invariante par translation. Comme $S^3$ est simplement connexe, d'apr\`es un r\'esultat 
  de~\cite{D} le groupe des isom\'etries $G$ est compact. Ceci implique que $I$ est compact et qu'il s'agit donc d'une g\'eom\'etrie riemannienne.
  
  Il reste \`a r\'egler le cas o\`u $G_{1}$ est $\widetilde{SL(2, \RR)}$. Consid\'erons $(X',Y,Z)$ une base de l'alg\`ebre de Lie $sl(2, \RR)$ avec les crochets de Lie usuels :
  $\lbrack X', Y\rbrack =Y$,  $\lbrack X', Z \rbrack= -Z$ et $\lbrack Z, Y \rbrack =2X'$. 
   
  Compl\'etons cette base en une base de $\mathcal G$ en ajoutant un g\'en\'erateur $T$
  de $\mathcal I$. 
  
  Consid\'erons l'action du groupe  d'isotropie $I$ en un point $x_{0}$ sur l'id\'eal $sl(2, \RR)$ qui s'identifie \`a $T_{x_{0}}M$. Cette action se fait par isomorphisme de l'alg\`ebre de 
  Lie $sl(2, \RR)$ et doit fixer un vecteur de $sl(2, \RR)$. Quitte \`a op\'erer un changement de base par un automorphisme int\'erieur de $sl(2, \RR)$ on peut supposer que
  le vecteur stabilis\'e par $I$ co\"{\i}ncide avec $X'$ (s'il est de norme non nulle pour la forme de Killing) ou bien avec $Y$ (s'il est de norme nulle pour la forme de Killing).
  
  Les seuls  isomorphismes de l'alg\`ebre de Lie $sl(2, \RR)$ qui fixent $X'$ \'etant de la forme $(X',Y,Z) \to (X', e^tY, e^{-t}Z)$, il vient que 
   l'action de $I$ (donc de $T$)  co\"{\i}ncide avec celle de $ad(X')$.  Le quatri\`eme champ de Killing $T$ est alors le sous-groupe \`a un param\`etre de 
   $\widetilde{SL(2, \RR)}$ engendr\'e par $X'$ et plong\'e diagonalement dans $\widetilde{SL(2, \RR)} \times    \widetilde{SL(2, \RR)}$.

  Nous venons de  d\'emontrer  que $G= \RR \times \widetilde{SL(2, \RR)}$ : la composante connexe de l'identit\'e du groupe des isom\'etries  est engendr\'ee par toutes les translations \`a gauche et par le sous-groupe \`a un param\`etre
  des  translations  \`a droite  engendr\'e par $X'$.
  
  De la m\^eme mani\`ere on constate que tous les isomorphismes de l'alg\`ebres de Lie $sl(2, \RR)$ qui fixent $Y$ sont engendr\'es par l'action adjointe de $Y$. Dans ce cas
   la composante connexe de l'identit\'e du groupe des isom\'etries  est engendr\'ee par toutes les translations \`a gauche et par le sous-groupe \`a un param\`etre
  des  translations  \`a droite  engendr\'e par $Y$.  
  
    Dans les deux  cas $G$ est un sous-groupe de $\widetilde {SL(2, \RR)}  \times \widetilde{SL(2, \RR)}$ qui est la composante connexe de l'identit\'e du groupe des isom\'etries de la m\'etrique lorentzienne \`a courbure sectionnelle constante n\'egative obtenue en consid\'erant la m\'etrique de Killing de $\widetilde{SL(2, \RR)}$. Notre g\'eom\'etrie n'est donc pas maximale : elle est subordonn\'ee \`a la g\'eom\'etrie
    de courbure sectionnelle constante n\'egative.\\
  
(2)   Supposons  que $\mathcal I$ est  contenu dans $\mathcal G$$_{1}$.  
  
  Ceci implique que les orbites de l'action de $\mathcal G$$_{1}$ sont de dimension $2$ et 
  que l'action de l'isotropie pr\'eserve une decomposition non triviale de l'espace tangent \`a $M$ (induite par la d\'ecomposition, $ad(I)$-invariante, $\mathcal G= \RR \oplus {\mathcal G}_{1}$). Il vient que l'isotropie $I$ est semi-simple, que le champ $X$ est engendr\'e par le facteur $\RR$ et que l'espace tangent aux orbites de ${\mathcal G}_{1}$ est $X^{\bot}$.
  
  Les orbites de $\mathcal   G$$_{1}$ sont alors  des surfaces lorentziennes homog\`enes et donc de courbure sectionnelle constante. Or, l'alg\`ebre de Lie de  $S^3$ n'est pas  l'alg\`ebre de Killing 
  d'une telle surface, tandis  que $sl(2, \RR)$  se r\'ealise bien comme l'alg\`ebre de Killing du plan de Sitter.

  Nous sommes alors dans le cas de la g\'eom\'etrie
  produit d\'ecrite au point (2), o\`u  l'isotropie est conjugu\'ee dans $\widetilde{SL(2, \RR)}$ \`a un sous-groupe \`a un param\`etre semi-simple.
  \end{demonstration}
     
 Nous passons \`a la preuve du corollaire~\ref{sl2}.

       \begin{demonstration}
  Par la proposition~\ref{semi-simple algebrique},  l'action de $G$ sur le mod\`ele $G/I$ pr\'eserve une m\'etrique de courbure sectionnelle constante n\'egative. Le lemme~\ref{completudes}
   et le r\'esultat de~\cite{Kli}  impliquent  que la $(G,G/I)$-structure de $M$ est compl\`ete. 
        
     Pour montrer que $M$ est g\'eod\'esiquement compl\`ete, il suffit alors de montrer que le mod\`ele $G/I$ est g\'eod\'esiquement complet. On utilisera les r\'esultats de~\cite{GL}
     sur la compl\'etude g\'eod\'esique des m\'etriques lorentziennes invariantes \`a gauche sur $\widetilde{SL(2, \RR)}$.
      
      Exprimons la forme quadratique $g$ par rapport \`a la forme de Killing $q$ de $sl(2, \RR)$ : $g(u,v)=q(\phi(u),v)$, o\`u $\phi$ est un endomorphisme $q$-sym\'etrique de
      $sl(2, \RR)$ et $u,v \in sl(2, \RR)$.
      De plus, $\phi$ est invariant par l'action d'un sous-groupe
      \`a un param\`etre semi-simple ou unipotent.
      
      Si le sous-groupe \`a un param\`etre est semi-simple, il vient que $\phi$ laisse stables les espaces propres du sous-groupe \`a un param\`etre semi-simple et il est
      donc diagonalisable. On v\'erifie que $\phi$ est $q$-sym\'etrique si et seulement si les deux valeurs propres qui correspondent aux deux vecteurs propres $q$-isotropes
      sont \'egales. Il vient que $\phi$ admet un espace propre de dimension au moins $2$ et on se trouve dans le cas (a) d'application du th\'eor\`eme $4$ de~\cite{GL}, qui assure 
      la compl\'etude g\'eod\'esique de $g$.
      
      Dans le cas restant,  le groupe \`a un param\`etre qui laisse invariant $g$ (et donc \'egalement $\phi$) est unipotent et $\phi$ s'exprime dans une base adapt\'ee
     (voir la preuve du lemme~\ref{invariance})  par une matrice de la forme $  \left(  \begin{array}{ccc}
                                                                 \lambda    &   \beta  & \alpha \\
                                                                 0     &  \lambda  &  - \beta  \\
                                                                 0     &   0  &  \lambda \\
                                                                 \end{array} \right),     \hbox{avec} \; \;   \alpha, \beta, \gamma \in \RR.$

                                                                  Cet
                                                                 endomorphisme est $q$-sym\'etrique si et seulement si $\beta =0$. Il vient que $\phi$ admet un espace propre de dimension au moins deux
                                                                 et on conclut comme avant.      
             \end{demonstration} 
   
   \section{Mod\`eles alg\`ebriques avec $G$  r\'esoluble}  \label{modeles algebriques 2}

  Dans cette section on classifie les g\'eom\'etries lorentziennes non-riemanniennes  $G/I$, avec $G$ r\'esoluble de dimension $4$ et $I$ de dimension $1$.

  Nous d\'eterminons d'abord $H$ : 
  
  \begin{proposition}      \label{algebrique}
  
  i) L'alg\`ebre d\'eriv\'ee $\lbrack \mathcal H, \mathcal H \rbrack$ est isomorphe \`a $\RR$.
  
  ii). Le groupe $H$ est isomorphe au groupe de Heisenberg ou bien au produit  $H=\RR \times AG$, o\`u $AG$ est le groupe affine.

 \end{proposition}

       Avant de passer \`a la preuve rappelons que  le groupe affine de la droite  $AG$ est un groupe de Lie  de dimension $2$
     qui peut \^etre vu comme l'ensemble des couples $(a,b) \in \RR^2$ muni de la multiplication $(a,b) \cdot (a',b')= (a+a', exp(a)b'+b).$

  \begin{demonstration}

  i)   Rappelons que l'alg\`ebre d\'eriv\'ee
  d'une alg\`ebre r\'esoluble est toujours nilpotente. Or,  $\mathcal H$ est suppos\'ee  r\'esoluble et donc $\lbrack \mathcal H, \mathcal H \rbrack$ est une alg\`ebre nilpotente de dimension strictement inf\'erieure  \`a $3$.
  
     L'alg\`ebre d\'eriv\'ee $\lbrack \mathcal H, \mathcal H \rbrack $ ne peut \^etre
     triviale. En effet, sinon $\mathcal H$ est ab\'elienne et donc l'action de l'isotropie $\mathcal I \subset \mathcal H$ sur $T_{x_{0}}F \simeq \mathcal H/ \mathcal I$ est triviale.
     Ceci est absurde car l'on a vu qu'un \'el\'ement de $\mathcal G$ agissant trivialement sur $F$ est trivial.
     
      Il vient que $\lbrack \mathcal H, \mathcal H \rbrack $  est isomorphe \`a $\RR$ ou \`a $\RR^2$.
     Supposons par l'absurde que l'alg\`ebre d\'eriv\'ee est isomorphe \`a $\RR^2$. 
     
     On d\'emontre d'abord que $\mathcal I$ est inclus dans $\lbrack \mathcal H, \mathcal H \rbrack$. Supposons par l'absurde  le contraire.
     Alors $\RR^2$ agit par isom\'etries sur chaque feuille $F$ du feuilletage $\mathcal F$.  Nous d\'emontrons alors que la restriction de la connexion $\nabla$ \`a $F$ est plate.
     Dans les coordonn\'ees locales $(x,h)$ de $F$ donn\'ees par l'action de $\RR^2$,  l'expression de $g$  est $dh^2$. Dans ces coordonn\'ees, $X$ est \'egalement invariant par translations et est de la forme $\frac{\partial}{\partial x}$, si $I$ est unipotent et de la forme $\frac{\partial}{\partial h}$, si $I$ est semi-simple.
     
     Un calcul simple montre que toute connexion sans torsion $\nabla$,  invariante par translations sur $\RR^2$,  et compatible avec $dh^2$ est donn\'ee par
     $\nabla_{\frac{\partial}{\partial h}}\frac{\partial}{\partial h}=a
\frac{\partial}{\partial x}$, $\nabla_{\frac{\partial}{\partial x}}\frac{\partial}{\partial x}=b
\frac{\partial}{\partial x}$ and $\nabla_{\frac{\partial}{\partial h}}\frac{\partial}{\partial x}=
\nabla_{\frac{\partial}{\partial x}}\frac{\partial}{\partial h}=c
\frac{\partial}{\partial x}$, avec    $a,b,c \in \RR$. L'invariance de $\nabla$ par $I$ implique qu'au moins deux des param\`etres $a,b$ et $c$ s'annulent. Ceci implique que
la courbure de $\nabla$ s'annule.

Par cons\'equent, le groupe des isom\'etries de ce mod\`ele est contenu dans le groupe affine de $\RR^2$. Il co\"{\i}ncide avec 
   $\RR \ltimes \RR^2$, o\`u l'action de l'isotropie $I \simeq \RR$ sur  $\RR^2$ est donn\'ee par 
$\left(  \begin{array}{cc}
                                                                 1   &   t\\
                                                                 0     &  1\\
                                                                 
                                                                 \end{array}  \right) $,  si $I$ est  unipotent, ou  par  $\left(  \begin{array}{cc}
                                                                 e^t   &   0\\
                                                                 0     &  1\\
                                                                 
                                                                 \end{array}  \right) $,  
                                                                si  $I$ est  semi-simple. Il vient que $H$ est isomorphe, ou bien \`a   Heisenberg,  ou bien \`a $AG \times \RR$. Dans les deux cas son groupe d\'eriv\'ee est de dimension un : absurde.

           Il vient donc que l'alg\`ebre d'isotropie $\mathcal I$
     est incluse  dans l'alg\`ebre d\'eriv\'ee $\lbrack \mathcal H, \mathcal H  \rbrack$.
     Les orbites de $\lbrack \mathcal H, \mathcal H \rbrack$ sur $F$ seront alors de dimension un. 
     
     Nous d\'emontrons que ces orbites sont les courbes int\'egrales sur $F$ de l'unique direction isotrope et que l'isotropie $I$ est unipotente. Pour cela, d\'esignons par $Y$ un g\'en\'erateur de l'alg\`ebre d'isotropie,   par $(Y,X')$ une base de l'alg\`ebre commutative $\lbrack \mathcal H, \mathcal H \rbrack$ et finalement
       par $(Y,X',Z)$ une base de $\mathcal H$. Comme $ad(Y) \cdot \mathcal H \subset \lbrack \mathcal H, \mathcal H \rbrack$ et que $ \lbrack \mathcal H, \mathcal H \rbrack$    est ab\'elienne, il vient que
     $ad(Y) \cdot X'=0$ et $ad(Y) \cdot Z=aX'$ (mod  $Y$). Donc,  l'action  de $ad(Y)$ sur l'espace tangent \`a $F$ est unipotente et admet $X'$ comme unique  direction (isotrope)  propre.  Or, $X'$ engendre  bien  l'espace tangent  aux orbites de $\lbrack \mathcal H, \mathcal H \rbrack$  sur $F$.
     
     D\'esignons par $\mathcal L$ l'alg\`ebre d\'eriv\'ee de $\mathcal G$. Nous avons les inclusions  $\mathcal I $$\subset \lbrack$$ \mathcal H$, $\mathcal H $\rbrack$  \subset
     \mathcal L$. La dimension de $\mathcal L$ est $2$ ou $3$ et les orbites de l'action de $\mathcal L$ sur le mod\`ele $G/I$ sont de dimension respectivement      $1$ ou $2$.
     
     Traitons d'abord le cas o\`u $\mathcal L$ est de dimension $3$ et ses orbites sont donc de dimension $2$. Comme l'unique $2$-plan de $T_{x_{0}}G/I$ pr\'eserv\'e par 
     l'isotropie unipotente est $X'^{\bot}$, il vient que les orbites de $\mathcal H$ et $\mathcal L$ co\"{\i}ncident.
     Ceci  est absurde car $\mathcal L$ est nilpotente (en tant qu'alg\`ebre d\'eriv\'ee
     d'une alg\`ebre r\'esoluble) et pas $\mathcal H$ (son alg\`ebre d\'eriv\'ee \'etant de dimension $2$).

     Il reste \`a traiter le cas o\`u $\mathcal L$ est de dimension $2$. Dans ce cas $\mathcal L =  \lbrack \mathcal H, \mathcal H \rbrack$. L'image de $\mathcal G$ par l'action de l'isotropie $ad(Y)$ est alors incluse dans $\mathcal L$. Cette  image est donc de dimension $2$ et de dimension $1$ dans $\mathcal G/ \mathcal I$ qui s'identifie \`a $T_{x_{0}}G/I$.
     Ceci est absurde car $ad(Y)$ doit \^etre de rang 2 (comme on l'a vu \`a la section~\ref{dynamique lorentzienne} en consid\'erant des bases adapt\'ees). On vient de montrer que  $\lbrack \mathcal H, \mathcal H \rbrack$ est de dimension $1$. \\
          
     ii). Consid\'erons un g\'en\'erateur $Z$ de l'alg\`ebre  $\lbrack \mathcal H, \mathcal H \rbrack$ et  son application adjointe $ad(Z) : \mathcal H$$ \to \RR Z$. Si l'application pr\'ec\'edente est identiquement nulle,
     autrement dit      si l'\'el\'ement $Z$ est central dans $\mathcal H$, alors $\mathcal H$ est nilpotente et $H$ est isomorphe au groupe de Heisenberg.
     
     Sinon, soit $X'$ un g\'en\'erateur du noyau de $ad(Z)$ et $Y$ un \'el\'ement de $\mathcal H$ tel que $Y,X',Z$ engendrent $\mathcal H$. Comme $Z$ n'est pas central, on peut
     supposer que $\lbrack Y,Z \rbrack =Z.$ Nous avons que $\lbrack X', Y \rbrack = \alpha Z$, pour $\alpha \in \RR.$ 
     
     Si $\alpha =0$, $H = \RR \times AG$, le centre
     \'etant engendr\'e par $X'$ et le deuxi\`eme facteur par $\RR Z   \oplus \RR Y$.
     
     Si $\alpha \neq 0$, il suffit de changer $X'$ en $X' + \alpha Z$ pour se ramener au cas pr\'ec\'edent.
     \end{demonstration}

  \subsection{Cas : $H=\RR \times AG$}
  
  \begin{proposition}     \label{classification algebrique}
   Si  $H= \RR \times AG$, l'isotropie $I$ est semi-simple  et engendr\'ee par le g\'en\'erateur des homoth\'eties dans $AG$. Le groupe $G$ est l'un des trois groupes suivants  :

\begin{enumerate}
\item
                                    $G= \RR \times SOL$,

                                    \item

                                     $G = \RR \ltimes Heis$, 
                                     
                                     \item 
                                     
                                     $G = \RR^2 \ltimes \RR^2$.
      \end{enumerate}                               
         
       Dans le cas $(2)$ l'action du premier facteur sur $Heis$,  dans une base $(X',Z,T)$ de $heis$   form\'ee par l'\'el\'ement central $X'$ et telle que $\lbrack T, Z \rbrack   =X'$,
       s'exprime par $(X',Z,T) \to (X', e^tZ,e^{-t}T)$.     
       
       Dans le cas $(3)$ l'action de la premi\`ere copie de $\RR^2$ sur la deuxi\`eme copie de $\RR^2$ est donn\'ee par les matrices       $\left(  \begin{array}{cc}
                                                                 e^t   &   0\\
                                                                 0     &  e^{-t}\\
                                                                 
                                                                 \end{array}  \right) $  et $\left(  \begin{array}{cc}
                                                                 1   &   0\\
                                                                 0     &  e^{-t}\\
                                                                 
                                                                 \end{array}  \right) $.
     \end{proposition}
     
     \begin{remarque}
  Comme le centre de $\mathcal G = \RR^2 \ltimes \RR^2$  est  trivial, le lemme~\ref{X} impliquera  que cette g\'eom\'etrie  ne se r\'ealise pas sur 
 des vari\'et\'es compactes.
\end{remarque}
     
     \begin{demonstration}
    Pour fixer les id\'ees supposons que $X'$ est un \'el\'ement central non trivial de $\mathcal H$
     et que $Y$ et $Z$ engendrent l'alg\`ebre de Lie du groupe affine : $\lbrack Y, Z \rbrack = Z$. Alors $X',Y,Z$ engendrent $\mathcal H$ et d\'esignons par $T$ un quatri\`eme
     g\'en\'erateur de $\mathcal G$.

 Nous montrons que,  quitte \`a appliquer un automorphisme de $\mathcal H$ qui envoie  $Y$ sur  $Y+aZ+bX'$, avec  $a,b \in \RR$, {\it l'alg\`ebre d'isotropie $\mathcal I$ est  $\RR Y$.}
 Pour cela, il suffit de montrer que $\mathcal I$ n'est pas incluse dans $\RR X' \oplus  \RR Z$.

Remarquons que  $ad(\alpha X' + \beta Z)(\mathcal H) \subset \RR X' \oplus  \RR Z$,   $\forall \alpha, \beta \in \RR$. Si $\mathcal I$ est $\RR (\alpha X' + \beta Z)$, alors 
l'action de  $ad(\alpha X' + \beta Z)$ sur  $T_{x_{0}}F \simeq \mathcal H /  \mathcal I$ est donn\'ee par une matrice de rang $1$. Ceci implique que $I$ n'est pas semi-simple.
Nous venons de prouver que dans le cas o\`u $I$ est semi-simple, alors $\mathcal I$ n'est pas incluse dans $\RR X' \oplus  \RR Z$ et on peut donc consid\'erer que $\mathcal I= \RR Y$.

D\'emontrons  la m\^eme chose dans le cas o\`u $I$ est unipotent. Observons d'abord que $\mathcal I \neq \RR X'$,  car l'\'el\'ement central $X'$ agit trivialement sur  $\mathcal H$
et donc aussi sur $\mathcal H / \mathcal I \simeq T_{x_{0}}F$.

Supposons par l'absurde que $\mathcal I \subset \RR X' \oplus \RR Z$. Quitte \`a appliquer un automorphisme de  $\mathcal H$ qui envoie  $Z$  sur  $Z + \alpha X'$, avec  $\alpha \in \RR$, supposons que 
 $\mathcal I = \RR Z$. Il vient que l'alg\`ebre commutative  $\RR X' \oplus \RR Y$ intersecte trivialement $\mathcal I$ et agit donc librement et transitivement sur $F$. Comme dans la preuve de la proposition~\ref{algebrique},  ceci implique que la  connexion induite sur $F$ par $\nabla$ est plate et que $H$ est le groupe de Heisenberg. Ceci est absurde : le groupe de Heisenberg et $\RR \times AG$ ne sont pas isomorphes.
                                                                 
  Il vient que, modulo un automorphisme de $\mathcal H$, on a $\mathcal I = \RR Y$. Mais ceci est impossible dans le cas unipotent. En effet, de nouveau $\RR X'  \oplus \RR Z$
     agit librement et transitivement sur $F$ et  $F$ est  plate. Comme avant,  $H$ est isomorphe au  groupe de Heisenberg, ce qui contredit notre hypoth\`ese.

     Nous venons de d\'emontrer que $H=\RR \times AG$ implique que l'isotropie est semi-simple et que $\mathcal I$ est engendr\'ee par $Y$.      
     
     Comme l'isotropie $ad(Y)$ fixe le vecteur $X'$ et dilate exponentiellement la direction engendr\'ee par $Z$, on peut choisir pour $T$ la deuxi\`eme direction isotrope
     du plan lorentzien $X'^{\bot}$. Alors $ad(Y) \cdot T= \lbrack Y, T \rbrack = -T + \alpha Y$, pour une certaine constante r\'eelle $\alpha$. Quitte \`a changer $T$ en $T - \alpha Y$,
     on  suppose que $\lbrack Y, T \rbrack = -T$.

     On montre que $\lbrack T,Z \rbrack =aX' + bY$, avec $a,b \in \RR$ et $\lbrack T, X' \rbrack=cT$.
     
    La relation de Jacobi  donne $\lbrack Y, \lbrack T, Z \rbrack \rbrack=  \lbrack  \lbrack Y, T \rbrack , Z \rbrack +  \lbrack T, \lbrack Y, Z \rbrack \rbrack =
     \lbrack -T, Z \rbrack + \lbrack T, Z \rbrack =0$. Ceci montre que $ \lbrack T, Z \rbrack $ centralise $Y$ et  doit donc appartenir \`a l'alg\`ebre engendr\'ee par $Y$ et $X'$.
     
     Par ailleurs, comme $X'$ et $Y$ commutent, $T$ est \'egalement un vecteur propre de $ad(X')$ (non seulement de $ad(Y)$). Ceci donne $\lbrack T, X' \rbrack= c T$, avec $c \in \RR$.
          
      Consid\'erons maintenant $\mathcal L = \lbrack \mathcal G, \mathcal G \rbrack.$ Comme     $\lbrack Y, Z \rbrack = Z$,        $\lbrack Y, T \rbrack = -T$ et        $\lbrack T,Z \rbrack =aX' + bY$, 
      nous avons  que $\mathcal L$ contient l'alg\`ebre de Lie engendr\'ee par $Z,T$ et $aX'+bY$.      On remarque que $\lbrack aX' + b Y, Z \rbrack = bZ$. Ceci implique $b=0$, car
      sinon l'alg\`ebre de Lie engendr\'ee par les vecteurs $aX'+ b Y$ et $Z$ est l'alg\`ebre du groupe affine qui ne peut pas \^etre contenue dans l'alg\`ebre nilpotente $\mathcal L$.
      
      Nous avons donc $b=0$ et $\lbrack T,Z \rbrack = aX'$. Comme $\lbrack T, aX' \rbrack =ac T$, le m\^eme argument implique que ou bien $a=0$, ou bien $c=0$. 
      
      Une renormalisation \'evidente
      permet de voir que $a \neq 0$ peut \^etre remplacer par $a=1$ et $c \neq 0$ peut \^etre remplacer  par $c=1$.
      
      Finalement on a trois alg\`ebres de Lie possibles pour $\mathcal G$.
      
      \begin{enumerate}
      
      \item 
      Si $a=0$ et $c=0$, les relations sont les suivantes : $\lbrack Y, Z \rbrack = Z, \lbrack Y, T \rbrack = -T, \lbrack T, Z \rbrack =0$ et $X'$ est central. Le groupe $G$ est alors
      $\RR \times SOL$,  o\`u $X'$ engendre le centre du groupe et l'alg\`ebre de Lie  de $sol$ est engendr\'ee par $Y,Z,T$. Le groupe $SOL$  peut \^etre vu comme un produit semi-directe $\RR \ltimes \RR^2$, o\`u l'action infinit\'esimale de $Y$ est donn\'ee dans
      la base $Z,T$, par la matrice $\left(  \begin{array}{cc}
                                                                 1   &   0\\
                                                                 0     &  -1\\
                                                                 
                                                                 \end{array}  \right) $. Le sous-groupe \`a un param\`etre $Y$ engendre l'isotropie qui se trouve enti\`erement incluse dans $SOL$.

      \item                                                                  
      Pour $a=1$ et $c=0$ on trouve  $\lbrack Y, Z \rbrack = Z, \lbrack Y, T \rbrack = -T, \lbrack T, Z \rbrack =X'$, $\lbrack T, X' \rbrack =0$ et le groupe correspondant est
      $\RR \ltimes Heis$ . 
      
      Le premier facteur $\RR$ est engendr\'e par le groupe d'isotropie $Y$.
      L'action du sous-groupe \`a un param\`etre engendr\'e par le premier facteur sur l'alg\`ebre de Lie de Heisenberg est $(X',Z,T) \to (X',e^t Z, e^{-t} T)$. L'\'el\'ement
      $X'$ qui engendre le centre de Heisenberg est donc \'egalement dans le centre de $G$.
      
      Le deuxi\`eme facteur, isomorphe au groupe de Heisenberg, intersecte trivialement le groupe d'isotropie et agit donc transitivement sur $G/I$.

      \item 
      Pour $a=0$ et $c=1$ on trouve : $\lbrack Y, Z \rbrack = Z, \lbrack Y, T \rbrack = -T, \lbrack T, Z \rbrack =0, \lbrack T, X' \rbrack =T$      et on a
      $G= \RR ^2 \ltimes \RR^2$. L'action adjointe du premier facteur, engendr\'e par $Y$ et $X'$ agit sur la copie de $\RR^2$ engendr\'ee par les \'el\'ements
      $Z$ et $T$ par les matrices $\left(  \begin{array}{cc}
                                                                 1   &   0\\
                                                                 0     &  -1\\
                                                                 
                                                                 \end{array}  \right) $      et respectivement $\left(  \begin{array}{cc}
                                                                 0   &   0\\
                                                                 0     &  -1\\
                                                                 
                                                                 \end{array}  \right) $.   
      \end{enumerate}  
      Ceci ach\`eve la preuve.
     \end{demonstration}
                                                                
  \subsection{Cas :  H isomorphe \`a Heisenberg}
  
  \begin{proposition}    \label{feuilles plates}
  
  i) L'isotropie $I$ est unipotente.
  
  ii) Les feuilles de $\mathcal F$ sont plates (pour la connexion de Levi-Civita) et $X$ est parall\`ele le long de $\mathcal F$.
  \end{proposition}
  
  \begin{demonstration} L'action de $I$ sur $\mathcal H/ \mathcal I \simeq T_{x_{0}}F$ ne pr\'eserve aucune d\'ecomposition non triviale de $T_{x_{0}}F$. Il en r\'esulte  que $I$ est unipotent.
  Aussi,  $I$ est diff\'erent du centre de Heisenberg car celui-ci  agit trivialement sur $\mathcal H /  \mathcal I$. Ceci implique qu'il existe des copies de $\RR^2$ transverse \`a $I$ dans $H$ et qui agissent donc librement et transitivement sur $F$. Comme dans la preuve de la proposition~\ref{algebrique},  ceci implique que $F$ est plate et que $X$ est parall\`ele en restriction \`a $F$.
  \end{demonstration}
  
  \begin{proposition}  \label{normal}
  $\mathcal H$ est     un id\'eal de $\mathcal G$. 
        \end{proposition}
        
   \begin{corollaire} Le feuilletage engendr\'e par l'action locale de $\mathcal H$ co\"{\i}ncide avec $\mathcal F$.
   \end{corollaire}
        
      \begin{demonstration}
        Consid\'erons $A$ un champ de Killing et $B$ un champ de vecteurs tangent \`a $X^{\bot}$.
        Nous montrons que $\lbrack A, B \rbrack =\nabla_{A}B - \nabla_{B}A$ est encore dans $X^{\bot}$. Le terme $\nabla_{A}B$ est dans $X^{\bot}$. En effet,
        $g(B,X)=0$ implique $g(\nabla_{A}B,X) =-g(\nabla_{A}X,B)=0$ (d'apr\`es la proposition~\ref{invariance},  $\nabla_{A}X$ est colin\'eaire \`a $X$).
        
        Par ailleurs, le champ de Killing $A$ pr\'eserve $X$ et donc $\nabla_{X}A = \nabla_{A}X$. Comme $\nabla_{\cdot} A$ est anti-sym\'etrique, on a $g(\nabla_{B}A, X)=-
        g(B, \nabla_{X}A)=-g(B, \nabla_{A}X)=0$, car $\nabla_{A}X$ est colin\'eaire \`a $X$. Le terme $\nabla_{B}A$ est donc \'egalement dans $X^{\bot}$, ce qui implique
        que  $\lbrack A, B \rbrack \in X^{\bot}$.
         \end{demonstration}
         
  {\bf Structure alg\`ebrique de $G$}.
  
 De ce qui pr\'ec\`ede r\'esulte  que $G$ est une extension de $H$. Consid\'erons, comme avant,  une base $\{X',Y,Z \}$ de $\mathcal H$ telle que $X'$ est central, $\lbrack Y, Z \rbrack =X'$ et $Y$
  engendre $\mathcal I$. Vu l'expression de l'action de $\mathcal I= \RR Y$ sur $\mathcal G/ \mathcal I \simeq T_{x_{0}}M$,  dans une base adapt\'ee, on peut choisir le quatri\`eme g\'en\'erateur $T$ de $\mathcal G$ tel que $ad(Y)T=-Z$ dans $\mathcal G/ \mathcal I$. Il vient que $\lbrack Y, T \rbrack = -Z + aY$, avec $a \in \RR$. Quitte \`a changer $Z$ en $Z-aY$,
  on obtient $\lbrack Y, T \rbrack =-Z$, tout en laissant les autres relations inchang\'ees.
  
  Comme l'action adjointe de $T$ sur $\mathcal H$ pr\'eserve le centre $\RR X'$ de $\mathcal H$, il existe une constante r\'eelle $c$ telle que $\lbrack T, X' \rbrack =c X'$.
  
  \begin{proposition}  \label{X' et  X}

i)  Il existe une fonction $H$-invariante $f$ sur $G/I$ telle que  $X'=fX$ ($f$ est localement d\'efinie sur $M$ et constante sur les feuilles de $\mathcal F$).

ii) $X$  est de Killing  (et  $f$ est  constante) si et seulement si  $c=0$.

iii) Dans la base  $\{ X', Z, Y \}$ de  $\mathcal H$ l'action de $T$ est donn\'ee par  $ad(T)=\left(  \begin{array}{ccc}
                                                                 c   &   m &    0\\
                                                                 0     &  c &    1\\
                                                                 0     &  n  &  0\\
                                                                 \end{array}  \right),$  avec  $m,n \in \RR$.
                                                                 
 (iv) Si  $c=0$ et  $n=0$, alors  $g$ est une m\'etrique plate,  invariante \`a  gauche sur $Heis$.
\end{proposition}

\begin{demonstration}

i) Le champ de droites  $\RR X'$ est $\mathcal G$-invariant. Or, dans le cas d'isotropie unipotente, le seul champ de droites $\mathcal G$-invariant dans $TM$ est celui engendr\'e par
$X$.  Donc,  $X'$ et $X$ sont en tout point colin\'eaires.

ii) Comme $X$ est $\mathcal G$-invariant, il est de Killing si et seulement s'il engendre le centre de $\mathcal G$.  Ainsi $X$ est de Killing si et seulement s'il est un multiple
constant de $X'$. Il est \'equivalent de dire que la fonction $f$ est constante ou  que $c=0$.

iii) On applique la relation de Jacobi aux \'el\'ements $Y,Z$ and $T$ de $\mathcal G$ et on constate que $ad(T)$ agit comme une d\'erivation sur $\mathcal G$ si et
seulement si $ad(T)Z= mX'+cZ+nY$, avec $m,n \in \RR$.

iv) Si $c=0$ et $n=0$, alors les \'el\'ements $X', Z$ et $T$ engendrent une sous-alg\`ebre de Lie de $\mathcal G$ isomorphe \`a l'alg\`ebre de Heisenberg (avec $X'$ comme
centre) qui agit librement et transitivement sur $M$. Il vient que $g$ est localement isom\'etrique \`a une m\'etrique invariante \`a gauche sur Heisenberg qui attribue la norme $0$
\`a l'\'el\'ement central. Or, ces m\'etriques lorentziennes sont plates \cite{Rah}.
\end{demonstration}

       \section{Classification : cas de l'isotropie  unipotente}   \label{isotropie unipotente}

 Ici, l'isotropie est suppos\'ee unipotente. D'apr\`es ce qui pr\'ec\`ede, nous sommes  dans
 le cas o\`u  $H$ est le groupe de Heisenberg. Rappelons que des constantes $c$ et $n$ ont \'et\'e definies \`a  la proposition~\ref{X' et  X}.  Le but de cette section est de d\'emontrer la
     
     \begin{proposition}  \label{unipotent}
     
     Soit  $(M,g)$ une vari\'et\'e lorentzienne compacte  localement model\'ee sur une g\'eom\'etrie  lorentzienne (non n\'ecessairement maximale) $G/I$,  avec $G$ r\'esoluble de dimension $4$  et $I$  sous-groupe \`a un param\`etre unipotent.
     
     (i) Si $c=0$, alors $g$ est localement isom\'etrique \`a une m\'etrique plate, invariante \`a gauche  sur $Heis$.
     
     (ii) Si $c \neq 0$, alors $M$ est localement model\'ee sur la g\'eom\'etrie Lorentz-SOL. Les seules r\'ealisations compactes de la g\'eom\'etrie  Lorentz-SOL sont, \`a rev\^etement fini pr\`es,  des quotients de $SOL$ par un r\'eseau $\Gamma$.
    \end{proposition}

   D\'emontrons d'abord la compl\'etude via les deux propositions suivantes : 
        
              \begin{proposition}  \label{feuilles completes}
              
              La $(H, H/I)$-structure est compl\`ete sur chaque feuille de $\mathcal F$. Les feuilles de $\mathcal F$ sont hom\'eomorphes  \`a des plans, \`a des cylindres ou \`a des tores.
              \end{proposition}
              
              \begin{demonstration}
              D\'emontrons  que la $(H, H/I)$-structure sur une feuille $F$ est compl\`ete.
        
         D\'esignons par $\widetilde F$ le rev\^etement universel de $F$ et par $\widetilde X$ l'image r\'eciproque du champ $X$ sur $\widetilde F$. L'application d\'eveloppante de notre
          $(H, H/I)$-structure envoie  $\widetilde F$ sur un ouvert de $\RR^2$ et le  champ de Killing $\widetilde X$ sur le  champ $ \frac{\partial}{\partial x}$ dans l'espace des
          param\`etres $(h,x)$ de $\RR^2$.
          
          Le champ $\widetilde X$ \'etant complet (car $X$ est d\'efini sur la vari\'et\'e compacte $M$), pour chaque orbite de $\widetilde X$  l'application d\'eveloppante r\'ealise un diff\'eomorphisme
          entre un ouvert connexe de $\widetilde F$    contenant l'orbite en question  et invariant par $\widetilde X$   et  une bande ouverte  horizontale 
          de la forme $\rbrack h - \epsilon , h + \epsilon \lbrack \times \RR$ dans $\RR^2$.         
         
                On peut construire un champ de vecteurs $V \in X^{\bot}$, globalement d\'efini  sur $M$, et de norme constante \'egale \`a $1$. Pour cela il suffit
          de choisir un champ de vecteurs dans $X^{\bot}$ qui n'est en aucun point collin\'eaire \`a $X$ (\`a rev\^etement double pr\`es, il n'y a pas de contrainte topologique pour cela) et de le diviser par sa norme (cette norme est non nulle car le noyau de la restriction de $g$ \`a $F$ est engendr\'e par $X$). Comme le feuilletage engendr\'e par $X$ est
          transversalement riemannien (pour la m\'etrique  induite par $g$ sur la transversale) et les orbites de $V$ sont param\'etr\'ees par la longueur, le flot de $V$ envoie n\'ecessairement une orbite de $X$ sur une autre orbite de $X$ (sans  respecter n\'ecessairement 
          le param\'etrage). L'image r\'eciproque $\widetilde V$ de $V$ par le rev\^etement universel de $F$ est un champ de vecteurs complet sur $\widetilde F$ (car $V$ complet sur $M$). Par connexit\'e, une orbite
          du flot de $\widetilde V$  intersecte chaque orbite de $\widetilde X$. Ceci est suffisant pour conclure \`a la compl\'etude de notre $(H, H/I)$-structure (pour plus de d\'etails, voir~\cite{Zeg} proposition 9.3).
          
          Passons \`a pr\'esent au type topologique des feuilles de $\mathcal F$. Comme les feuilles admettent le champ de vecteurs non singulier $X$, si elles sont compactes
          elles sont n\'ecessairement hom\'eomorphes \`a des tores. Si une feuille $F$ est ouverte, son  groupe fondamental $\pi_{1}(F)$ agit sur le rev\^etement universel $\RR^2$
          en pr\'eservant les flots des champs de vecteurs $\widetilde X$ et $\widetilde V$. Il vient que $\pi_{1}(F)$ agit sur la transversale du feuilletage engendr\'e par $\widetilde X$ en
          pr\'eservant le flot de $\widetilde V$ : ceci donne un morphisme de $\pi_{1}(F)$ dans $\RR$. Le noyau de ce morphisme est form\'e par les \'el\'ements de $\pi_{1}(F)$
          qui fixent les feuilles donn\'ees par l'action de $\widetilde X$ : ce sous-groupe de $\pi_{1}(F)$ commute avec le flot de $\widetilde X$ et est donc \'egalement commutatif. Il vient
          que $\pi_{1}(F)$ ne peut \^etre un groupe libre, sans \^etre isomorphe \`a $\ZZ$ ou trivial. Dans ce cas $F$ est n\'ecessairement  hom\'eomorphe \`a un cylinde ou \`a un plan.
     \end{demonstration}
           
\begin{proposition}

La $(G,G/I)$-structure de $M$ est compl\`ete.
\end{proposition}

\begin{demonstration}
La $(G,G/I)$-structure de $M$ est une combinaison de  la $(H,H/I)$-structure des feuilles de  $\mathcal F$ et de la structure transverse de translation du feuilletage $\mathcal F$, qui est model\'ee sur  $G/H \simeq {\bf R}$.
Par compacit\'e de $M$, la structure transverse  de $\mathcal F$  est compl\`ete~\cite{Mo}. Comme les feuilles de $\mathcal F$ sont compl\`etes par la proposition~\ref{feuilles completes}, il vient que la $(G,G/I)$-structure de $M$ est \'egalement compl\`ete.
\end{demonstration}

\begin{lemme}  \label{holonomie}

i) Le groupe d'holonomie $\Gamma$ n'est pas contenu dans $H$.

ii) $\Gamma$ n'est pas ab\'elien.

iii) Si $c=0$ et $\Gamma$ est  nilpotent, alors $g$ est plate,  invariante \`a gauche sur $Heis$.
\end{lemme}

\begin{demonstration} Consid\'erons l'adh\'erence de Zariski $\overline{\Gamma}$ de $\Gamma$ dans $G$. Il s'agit d'un sous-groupe alg\'ebrique de $\Gamma$, qui admet donc un nombre fini de
composantes connexes. Quitte a consid\'erer un rev\^etement fini de $M$, on peut supposer que $\overline{\Gamma}$  est connexe. Par compl\'etude, l'application
d\'eveloppante de la $(G, G/I)$-structure de $M$ fournit une application surjective de $M$ dans le double quotient $\Gamma \backslash G/I$, et donc \'egalement une application
de $M$ sur $\overline{\Gamma} \backslash G/I$.

i) Supposons par l'absurde que $\overline{\Gamma} \subset H$. Nous avons alors une application de $M$ sur  le quotient de  $H \backslash G/I$. Or, ce quotient, qui s'identifie \`a la transversale $\RR$ du feuilletage $\mathcal F$,
 est un espace s\'epar\'e et non compact qui ne peut \^etre l'image continue du compact $M$ : absurde.
 
 ii) Supposons par l'absurde que  $\Gamma$ et, par cons\'equent, aussi $\overline{\Gamma}$ est ab\'elien. Un calcul imm\'ediat montre que le centralisateur  dans $\mathcal H$ d'un \'el\'ement appartenant \`a la diff\'erence  $\mathcal G \setminus  \mathcal H$ est de dimension \'egale \`a $1$, engendr\'e par $X'$, si $c=0$,  ou par un \'el\'ement de $\mathcal H$ non contenu dans $\RR   X' \oplus \RR   Y$, si $c \neq 0$. Il vient que
 $\overline{\Gamma}$ est de dimension au plus $2$.
 
 Supposons d'abord que $\overline{\Gamma}$ est un sous-groupe \`a un param\`etre (non contenu dans $H$). Dans ce cas le double quotient $\overline{\Gamma} \backslash G/I$
 s'identifie \`a $H/I$ qui n'est pas compact : absurde.
      
      Supposons maintenant  que $\overline{\Gamma}$
      est de dimension $2$. 
      
      Si $c \neq 0$,  alors $\overline{\Gamma}$  intersecte $H$ selon un sous-groupe \`a un param\`etre engendr\'e par un \'el\'ement $Z'$ de $\mathcal H$ non contenu
      dans $\RR   X' \oplus \RR  Y$. Le groupe $\overline{\Gamma}$ est engendr\'e alors par le sous-groupe \`a un param\`etre associ\'e \`a $Z'$ et par n'importe quel sous-groupe \`a un param\`etre $l$,  contenu
      dans $\overline{\Gamma}$ et transverse \`a $H$. Le quotient $\overline{\Gamma} \backslash G/I$ s'identifie alors au quotient \`a gauche de $H/I$ par
      le sous-groupe \`a un param\`etre engendr\'e par $Z'$. 
      
      Comme $Z^\prime \notin \RR   X' \oplus \RR   Y$, l'alg\`ebre commutative $\RR X' \oplus \RR Z'$ agit transitivement sur les feuilles de $\mathcal F$. 
      En particulier, le champ de Killing $Z'$ est parall\`ele,   et en restriction \`a la feuille $H/I$ de $\mathcal F$, une orbite du champ de vecteurs
      complet $X$ constitue une transversale totale au feuilletage d\'efini par $Z'$. Ceci implique que le quotient de $H/I \simeq \RR^2$ par l'action de $Z'$ s'identifie \`a la transversale
      $\RR$ du feuilletage trivial d\'efinit par $Z'$. Or, $\RR$ ne peut pas \^etre l'image continue du compact $M$ : absurde.
      
      Si $c=0$, $\overline{\Gamma}$ est engendr\'e  par $X'$ et par n'importe quel sous-groupe \`a un param\`etre $l$ contenu
      dans $\overline{\Gamma}$ et transverse \`a $H$. Le quotient $\overline{\Gamma} \backslash G/I$ s'identifie alors au quotient \`a gauche de $H/I$ par
      le sous-groupe \`a un param\`etre engendr\'e par le centre $X'$. Ce quotient est isomorphe \`a $\RR$ qui n'est pas compact : absurde.
      
  iii) Comme $G$ n'est pas nilpotent, il vient que $\overline{\Gamma}$ est de dimension $3$, isomorphe
  au groupe de Heisenberg.
  On a donc que l'alg\`ebre de Lie d\'eriv\'ee de l'alg\`ebre de Lie de $\overline{\Gamma}$     est $\RR  X'$. Or, si $T+aY+bZ$ et $T+a'Y+b'Z$ sont deux \'el\'ements lin\'eairement ind\'ependants du quotient de l'alg\`ebre de Lie de  $\overline{\Gamma}$ par $\RR  X'$ (avec $a,b$ et $a',b'$ des constantes r\'eelles), leur crochet qui vaut $(a-a')Z    + (b-b') nY$ n'est  nul que si $n=0$.  Donc $n=0$ et la proposition~\ref{X' et X} implique alors que $g$ est plate.
  \end{demonstration}
  
  \begin{remarque} Dans la preuve pr\'ec\'edente  nous consid\'erons l'adh\'erence de 
   Zariski $\overline{\Gamma}$ de $\Gamma$ dans $G$. Cette  op\'eration est l\'egitime  si $G$ est un groupe de Lie alg\'ebrique r\'eel. Dans  les exemples qu'on \'etudie
   $G$ n'est pas toujours un groupe alg\'ebrique r\'eel, mais il est isomorphe en tant que groupe de Lie \`a la composante connexe de l'identit\'e d'un groupe alg\'ebrique r\'eel.
  \end{remarque}
  
  Passons maintenant \`a la preuve de la proposition~\ref{unipotent} :
        
\begin{demonstration}

(i) Supposons d'abord $c=0$.
 D'apr\`es la proposition~\ref{X' et X}, $X$ est alors un champ de Killing isotrope  globalement d\'efini.

     Il est prouv\'e dans~\cite{Zeg} (section 14) que le champ de Killing  $X$ est n\'ecessairement {\it \'equicontinu}. Par d\'efinition, ceci signifie que  l'adh\'erence du  flot de $X$ dans le groupe des hom\'eomorphismes
     de $M$ est un groupe compact. 
     
     Par cons\'equent, le flot de $X$ pr\'eserve \'egalement une m\'etrique riemannienne sur $M$. Nous pouvons donc utiliser la classification 
    des champs de Killing riemanniens sur les  vari\'et\'es compactes  de dimension $3$~\cite{Ca}. Dans notre cas,  le feuilletage engendr\'e par $X$ est transversalement de Lie, 
    localement model\'e (transversalement) sur le quotient de $G$ par son centre. Un r\'esultat  de~\cite{Mo} (th\'eor\`eme 4.2) affirme alors que les adh\'erences des orbites de $X$ ont toutes la
    m\^eme dimension. Selon la dimension de l'adh\'erence des orbites de $X$ les   situations possibles sont les
    suivantes~\cite{Ca} :\\

    1) Si les orbites de $X$ sont denses dans $M$, alors $M$ s'identifie \`a un tore $T^3$  sur lequel le flot de $X$ est \`a orbites denses et le feuilletage engendr\'e par $X$ est lin\'eaire. Comme le groupe fondamental de $T^3$ est ab\'elien, il vient que $\Gamma$ est ab\'elien, ce qui est impossible par le lemme~\ref{holonomie}.\\

  2) Si les adh\'erences des orbites sont de dimension $2$, il est montr\'e \'egalement dans~\cite{Ca} (voir le th\'eor\`eme 1 de la section  III.A et le corollaire 4 de la section III.B) que  $M$ est un tore $T^3$ de dimension $3$.       On conclut comme pr\'ec\'edemment .\\

     3)   Il reste \`a traiter le cas o\`u les orbites de $X$ sont ferm\'ees. Dans ce cas, \`a rev\^etement fini pr\`es,  $M$ est une fibration principale sur un tore $T^2$ ou bien sur une sph\`ere $S^2$ avec le champ de Killing  $X$ qui engendre la fibration principale~\cite{Ca} (le cas des fibr\'es de Seifert non triviaux est \'elimin\'e gr\^ace \`a l'existence de
     la structure transverse : l'holonomie d'une orbite de $X$ ne peut \^etre une rotation rationnelle non triviale).
     
         Comme les orbites de $X$ sont ferm\'ees,  l'holonomie $\Gamma$ intersecte non trivialement
      le centre de $H$ et cette intersection est un sous-groupe discret  isomorphe \`a $\ZZ$ (image du groupe fondamental de la fibre par le morphisme d'holonomie).
      
     Ainsi l'holonomie $\Gamma$ est une extension centrale d'un groupe ab\'elien (l'image du groupe fondamental de $T^2$) et est donc nilpotente. Le lemme~\ref{holonomie}
        implique alors que $g$ est plate, invariante \`a gauche sur $Heis$.\\

     (ii) Consid\'erons maintenant le cas $c \neq 0$.

  Prouvons d'abord que    $\Gamma \cap H$ est un sous-groupe $\Delta$ non trivial de $\Gamma$.
  
 Supposons par l'absurde que $\Gamma \cap H =\{1\}$.
       Il vient donc que $\Gamma$ s'injecte dans $G/H \simeq \RR$, ce qui  implique que $\Gamma$ est commutatif, et  contredit    le lemme~\ref{holonomie}.

      Comme les feuilles de $\mathcal  F$ sont compl\`etes, une
     telle feuille s'identifie au quotient de $H/I$ par l'action  (\`a gauche) de $\Delta$. Il vient que $\Delta$ s'identifie au groupe fondamental d'une feuille de $\mathcal F$ et est un sous-groupe discret de $H$ isomorphe \`a $\ZZ$ (les feuilles sont  hom\'eomorphes 
     \`a des cylindres) ou \`a
     $\ZZ \oplus \ZZ$ (les feuilles sont hom\'eomorphes
       \`a des tores).

      Consid\'erons $\gamma$ un \'el\'ement de $\Gamma$
  ayant une projection non-triviale sur $G/H$. Pour ce qui va suivre, on ne fera plus appel aux crochets exacts de $T$ avec $Y$ et $Z$. On se permettra donc de modifier  $T$ 
   en lui ajoutant un \'el\'ement de ${\mathcal H}$. On peut ainsi supposer que 
  $\gamma = \exp(\alpha T) $,  pour un certain  $\alpha \in {\bf R}$. Quitte \`a changer $T$ par $-T$ et $\gamma$ par $\gamma^{-1}$, on  suppose que   $\alpha >0$ et $c <0$.

Consid\'erons maintenant   l'action de $\gamma$
  sur $H$, et en particulier sur $\Delta$ et sur sa fermeture de Zariski 
  $\overline{\Delta}$ dans $H$. \\

  Consid\'erons d'abord le cas $\Delta \simeq \ZZ$. Comme ${\bf Z}$ n'admet pas d'automorphismes autres que $z \to -z$, \`a indice 2 pr\`es, 
  $\Gamma$ agit trivialement sur $\Delta$, et donc \'egalement sur le groupe \`a un param\`etre $\overline{\Delta}$. Comme $c \neq 0$, l'\'el\'ement $\gamma$ ne
  pr\'eserve pas $X'$ et donc $\overline{\Delta}$ admet un g\'en\'erateur infinit\'esimal $Z' \in \mathcal H$, non colin\'eaire \`a $X'$.

  Il vient aussi que $\Gamma$ est contenu dans le centralisateur $C$ de $\overline{\Delta}$ dans $G$. Le sous-groupe $C$ de $G$ admet un centre de dimension au moins $1$ et  est donc de dimension au plus $3$ (comme $c\neq 0$, le centre de
  $G$ est trivial). 
  
  L'alg\`ebre de Lie $\mathcal C$ 
  de $C$ contient $\RR X' \oplus \RR Z'$, et  doit \'egalement  contenir $T$.
  Ainsi ${\mathcal C}$ est d\'efinie par les relations 
  $\lbrack T, Z' \rbrack = \lbrack X', Z' \rbrack  = 0$ et $\lbrack T, X' \rbrack  = cX'$.  Comme  $\Delta$ et (donc)  $\overline{\Delta}$ agissent proprement sur $H/I$, il s'en suit que $\RR Z'$
  intersecte trivialement l'isotropie $\RR Y$. Il en d\'ecoule que le groupe
  $C$ agit librement et transitivement sur le mod\`ele   $G/I$ et que  $M$ admet une 
  $(C,C)$-structure. En particulier, $\Gamma$ est un r\'eseau (cocompact) de $C$ et $M$ s'identifie au quotient $\Gamma \backslash C$. Or, $C$ est isomorphe \`a ${\bf R} \times AG$ et n'admet pas de r\'eseau car il n'est pas unimodulaire.  Ceci est absurde.\\

Consid\'erons \`a pr\'esent le cas $\Delta \simeq  \ZZ^2$. 
 Alors $\overline{\Delta}$ est ab\'elien  et de dimension $2$. Il doit contenir  le centre de $H$. Notons son alg\`ebre de Lie $\RR X' \oplus \RR  Z'$, o\`u $Z'$ est un \'el\'ement de $\mathcal H$  non colin\'eaire \`a $X'$.

Les feuilles de $\mathcal F$ sont des tores.  La th\'eorie g\'en\'erale des structures transverses des feuilletages~\cite{Mo} assure que 
 $M$ est un fibr\'e en tores sur un cercle et l'image de $\Gamma$ par la projection $G \to G/H$ est un sous-groupe discret isomorphe \`a $\ZZ$. On peut  alors supposer    que l'\'el\'ement choisi $\gamma$ engendre la projection de $\Gamma$
sur $G/H$. Il vient alors que $\Gamma$ est inclus dans le groupe engendr\'e par $\gamma$ et $\Delta$.

 Rappelons que 
$\gamma$ agit sur le plan $\overline{\Delta}$ en pr\'eservant le r\'eseau 
$\Delta$.  Mais  le groupe pr\'eservant un r\'eseau est unimodulaire. Il en d\'ecoule que $\gamma$ agit sur $\overline{\Delta}$ en pr\'eservant le volume. Comme $\lbrack T, X' \rbrack = cX'$, l'action de $\gamma$ admet 
une valeur propre \'egale \`a $e^{\alpha c}<1$. Il s'ensuit que l'action de 
$\gamma$   sur  $\RR X' \oplus \RR  Z'$
 est diagonalisable, avec deux valeurs propres $e^{\alpha}c$ et $e^{-\alpha c}$.  Ainsi, l'alg\`ebre de Lie engendr\'ee par 
 $X',Z'$  et $T$ est isomorphe \`a  $sol$. Comme avant, $\RR X' \oplus \RR Z'$   intersecte trivialement l'isotropie et  $SOL$ agit librement et transitivement sur $G/I$. Par la proposition~\ref{Lorentz.SOL1},  {\it nous  sommes en pr\'esence d'une 
 g\'eom\'etrie Lorentz-SOL},  car l'alg\`ebre  d\'eriv\'ee engendr\'ee par $X'$ et $Z'$ est d\'eg\'en\'er\'ee et la direction propre $\RR X'$ de $ad(T)$ est isotrope. Par ailleurs, nous avons
 montr\'e que l'holonomie $\Gamma$ est contenue dans $SOL$. Il vient que $M$ poss\`ede une $(SOL,SOL)$-structure, ce qui implique que $M$ est un quotient de $SOL$ par un r\'eseau cocompact.
 \end{demonstration}

  \section{Classification : cas de l'isotropie semi-simple}   \label{cas semi-simple}
  
  L'isotropie est suppos\'ee semi-simple. Traitons d'abord le cas
  
  \subsection{$G$ r\'esoluble}
  
  Le but de cette partie  est de d\'emontrer la 
  
      \begin{proposition}  \label{semi-simple}
     
     Soit  $(M,g)$ une vari\'et\'e lorentzienne compacte  localement model\'ee sur une g\'eom\'etrie  lorentzienne (non n\'ecessairement maximale) $G/I$,  avec $G$ r\'esoluble de dimension $4$  et $I$  sous-groupe \`a un param\`etre semi-simple. Alors $g$ est localement isom\'etrique, ou bien \`a une m\'etrique plate, invariante \`a gauche sur $SOL$, ou bien \`a la m\'etrique Lorentz-Heisenberg.
     
      Les seules r\'ealisations compactes de la g\'eom\'etrie  Lorentz-Heisenberg  sont, \`a rev\^etement fini  pr\`es, des quotients de $Heis$ par un r\'eseau $\Gamma$.
      \end{proposition}

    On utilise la classification des alg\`ebres $\mathcal G$ obtenue dans la proposition~\ref{classification algebrique}.

   \subsubsection*{Cas $G= \RR \times SOL$}
   
   \begin{proposition}   \label{non maximale}
   
   La g\'eom\'etrie $(G,G/I)$ repr\'esente une m\'etrique invariante \`a gauche sur $SOL$. Elle n'est donc pas maximale.
   \end{proposition}
   
   \begin{demonstration}
      Dans ce cas $G= \RR \times SOL$, o\`u $SOL$ est engendr\'e par   $\{ Z,T,Y \}$ (avec les relations $\lbrack Y, Z \rbrack = Z, \lbrack Y, T \rbrack = -T$ et $\lbrack T, Z \rbrack =0$) et le centre est engendr\'e par $X'$. L'alg\`ebre de Lie engendr\'ee par
  $\{ X',Z,T  \}$ est ab\'elienne et agit librement transitivement sur $G/I$. La m\'etrique $g$  est donc plate. Cette m\'etrique  s'identifie \`a une m\'etrique lorentzienne invariante \`a gauche  sur le groupe $SOL$ engendr\'e par les \'el\'ements $\{Y+X',Z,T \}$.
  
   Le mod\`ele  $G/I$ ne repr\'esente pas une g\'eom\'etrie lorentzienne maximale : la  g\'eom\'etrie maximale correspondante est  la g\'eom\'etrie Minkowski. 
  \end{demonstration}

   \subsubsection*{ Cas $G= \RR \ltimes Heis$}
          
          \begin{proposition} La g\'eom\'etrie $(G,G/I)$ est la g\'eom\'etrie Lorentz-Heisenberg.
          \end{proposition}
          
          \begin{demonstration}
             Le centre de $G$ est non trivial et engendr\'e par un \'el\'ement central $X'$ de $Heis$.  Ceci donne l'existence d'un champ de Killing globalement d\'efini sur $M$ de norme constante  positive  et qui est pr\'eserv\'e par l'action de $\mathcal G$.  Il co\"{\i}ncide donc avec un multiple du  champ de Killing  $X$ stabilis\'e par l'isotropie.
           
          Le  deuxi\`eme facteur $Heis$ agit librement et transitivement sur $G/I$. Ceci implique que la m\'etrique $g$ s'identifie localement \`a une  m\'etrique invariante \`a gauche sur $Heis$ qui attribue  \`a l'\'el\'ement central  $X'$ une norme positive.    La vari\'et\'e $M$ est alors localement model\'ee sur la g\'eom\'etrie Lorentz-Heisenberg (voir section \ref{metriques invariantes}).
           \end{demonstration}

 D\'emontrons maintenant la compl\'etude et la rigidit\'e de Bieberbach des r\'ealisations compactes de la  g\'eom\'etrie Lorentz-Heisenberg.
    
    \begin{demonstration}
     Dans le cas o\`u le champ de Killing $X$ (de norme \'egale \`a $1$) est non-\'equicontinu, il a \'et\'e montr\'e dans~\cite{Zeg} que $X$ est n\'ecessairement un flot d'Anosov dont les feuilletages
    stables et instables sont les deux droites isotropes de $X^{\bot}$. Comme $\mathcal G$ est r\'esoluble,  $M$ est n\'ecessairement une suspension d'un diff\'eomorphisme hyperbolique d'un tore et $g$ est plate (voir~\cite{Zeg}). Ceci est absurde car les seules m\'etriques lorentziennes plates et invariantes par translations sur $Heis$
    attribuent \`a l'\'el\'ement central $X'$ la norme $0$~\cite{Rah}.
    
    Il reste \`a   analyser le cas o\`u le flot de $X$ est \'equicontinu. Le feuilletage engendr\'e par $X$ admet une structure transverse qui
    est \`a la fois lorentzienne et riemannienne.  En particulier, le feuilletage engendr\'e par $X$ est transversalement de Lie et  la dimension de l'adh\'erence
    d'une orbite de $X$ ne d\'epend pas de l'orbite choisie~\cite{BMT,Mo}. Selon la valeur de cette dimension nous avons les possibilit\'es  suivantes~\cite{Ca} (voir \'egalement le th\'eor\`eme 4.2 dans~\cite{BMT}) :\\
    
   1)  Si les orbites de $X$ sont denses dans $M$ alors $M$ s'identifie \`a un tore $T^3$  sur lequel le flot de $X$ est \`a orbites denses et le feuilletage engendr\'e par $X$ est lin\'eaire. 
    L'adh\'erence du flot de $X$ dans le groupe des hom\'emorphismes de $M$ est alors un groupe de Lie ab\'elien compact qui agit transitivement par isom\'etries, aussi bien pour la m\'etrique
    lorentzienne $g$, que pour une m\'etrique riemannienne. Ceci implique que l'adh\'erence du flot de $X$ est un tore $T^3$ qui agit simplement transitivement et par isom\'etries sur $M$. Il en r\'esulte  que $g$ est \'egalement plate~: absurde.\\
    
    2) Si les adh\'erences des orbites de $X$ sont de dimension $2$,  il est montr\'e dans~\cite{Ca} (voir le th\'eor\`eme 1 de la section  III.A et le corollaire 4 de la section III.B)
    que $M$ est un tore $T^3$. Par  ailleurs, l'adh\'erence du flot de $X$ dans les hom\'eomorphismes de $M$ est un tore $T^2$ qui agit par isom\'etries pour la
    m\'etrique lorentzienne $g$, ce qui entra\^{\i}ne  l'existence de deux champs de Killing p\'eriodiques et qui commutent  sur $M$.  Ceci implique que le groupe de holonomie
    $\Gamma$ contient un sous-groupe discret isomorphe \`a $\ZZ^2$.
    
    Comme le groupe fondamental de $T^3$ est ab\'elien, il vient que l'holonomie $\Gamma$ est \'egalement ab\'elienne.
   Comme $\Gamma$ contient un sous-groupe discret isomorphe \`a $\ZZ^2$,  son adh\'erence  de Zariski $\overline{\Gamma}$ est un sous-groupe de Lie
    ab\'elien de dimension au  moins $2$ de $G$. Quitte \`a consid\'erer un rev\^etement fini  de $M$, $\overline{\Gamma}$ sera suppos\'e
    connexe.\\
    
     Supposons par l'absurde que $\Gamma$ n'est pas contenu dans $Heis$.
    
    On constate que le centralisateur dans l'alg\`ebre de Lie $heis$ d'un \'el\'ement de $\mathcal G$ qui ne se trouve pas dans $heis$
    est exactement $\RR   X'$.     Donc $\overline{\Gamma}$  est de dimension $2$. Plus pr\'ecis\'ement, $\overline{\Gamma}$ est une copie de  $\RR^2 \subset G$ engendr\'ee par $\RR   X'$  et par n'importe quel  sous-groupe
    \`a un param\`etre $l$ contenu dans $\RR^2$ et transverse \`a $\RR   X'$.

       L'image de  l'application d\'eveloppante est  un ouvert $U$ de $G/I$. Nous avons donc une application surjective de $M$ dans le quotient de $U$  par l'action de $\overline{\Gamma}$. Comme le champ de Killing $X$ est globalement d\'efini sur $M$, l'ouvert $U$ est invariant par 
    l'action de $X'$.

    Pour comprendre le quotient \`a gauche de $G/I$ par $\overline{\Gamma}$, il est ais\'e de regarder d'abord le quotient de $G$ par son centre, engendr\'e par  $\RR   X'$. On constate
    que ce quotient est isomorphe au groupe $SOL$ et que le quotient \`a gauche  de $G/I$ par $\overline{\Gamma}$ s'identifie au quotient \`a gauche de $SOL/I$ par l'image
    $l'$ de $l$ dans $G/I \simeq SOL$. Le mod\`ele $SOL/I$ est le plan de Minkowski     et l'image de l'ouvert satur\'e $U$ dans $SOL/I$ est un ouvert $U'$.
    
    L'application d\'eveloppante fournit alors une application surjective de $M$ dans le quotient de $U'$ par le champ de Killing $l' \subset SOL$, non contenu dans les translations
    (car $l$ non contenu dans $Heis$).
    
    La contradiction recherch\'ee viendra du fait que la norme de $l'$ (constante sur les orbites de $l'$) descend en une fonction continue sur $M$ sans extremas locaux.
    
    En effet, pour pr\'eciser cette id\'ee consid\'erons des coordonn\'ees $(z,t)$ sur $SOL/I$ dans lesquelles la m\'etrique lorentzienne s'exprime $q= dz   dt$.
     Comme  $l'$ est un champ de Killing qui n'est pas une translation pure, alors quitte \`a le multiplier par une constante, $l'=( z \frac{\partial}{\partial z}  - t \frac{\partial}{\partial t})
    + a \frac{\partial}{\partial z} + b \frac{\partial}{\partial t}$, avec $a,b \in \RR$. Dans ce cas l'expression de $q(l')$ est la fonction $(z+a)(b-t)$ qui  n'admet pas d'extrema local sur $\RR^2$. 
    
    Pourtant la fonction $q(l')$ descend bien en une fonction continue non constante sur la vari\'et\'e compacte $M$ : elle
    devrait donc admettre sur l'ouvert $U'$ au moins un minimum  et un maximum. Cette contradiction ach\`eve la preuve de $\Gamma \subset Heis$ (\`a indice fini  pr\`es), dans le cas o\`u l'adh\'erence
    d'une orbite de $X$ est de dimension $2$.\\

    3)  Il reste \`a r\'egler le cas o\`u toutes les orbites de $X$ sont ferm\'ees, de dimension $1$. Dans ce cas $M$ est  une fibration principale en cercles sur un tore $T^2$, la fibration principale \'etant engendr\'ee par $X$. Les orbites de $X$ \'etant ferm\'ees,  le groupe d'holonomie $\Gamma$ intersecte 
    le sous-groupe \`a un param\`etre engendr\'e par l'\'el\'ement central $X'$  selon un sous-groupe discret isomorphe \`a $\ZZ$. Il vient que le sous-groupe \`a un param\`etre
    engendr\'e par $X'$ est contenu dans $\overline{\Gamma} \cap Heis$, o\`u $\overline{\Gamma}$ est l'adh\'erence de Zariski de $\Gamma$.
    
   Il en r\'esulte que l'holonomie $\Gamma$ est une extension centrale d'un groupe ab\'elien (l'image du groupe fondamental de $T^2$ par le morphisme d'holonomie). Le groupe
   $\Gamma$ est donc nilpotent.\\
    
    Supposons par l'absurde que $\Gamma$ n'est pas contenu dans $Heis$.  Alors l'adh\'erence de Zariski $\overline{\Gamma}$ de $\Gamma$ est un sous-groupe nilpotent de 
    $G$ qui contient le centre $\RR   X'$ et qui n'est pas contenu dans $G$. Comme avant, on peut supposer $\overline{\Gamma}$ connexe.
    
    On montre que  $\overline{\Gamma}$ est n\'ecessairement de dimension $2$. D'abord, la dimension ne peut \^etre $4$ car $G$ n'est pas nilpotent. 
    
    Le centralisateur dans l'alg\`ebre de Lie $heis$ d'un \'el\'ement de $\mathcal G$ qui ne se trouve pas dans $heis$
    est exactement $\RR   X'$. Par ailleurs, si l'alg\`ebre de Lie de $\overline{\Gamma}$ est suppos\'ee de dimension $3$, alors l'intersection de l'alg\`ebre de Lie
    de $\overline{\Gamma}$ avec $heis$ est de dimension $2$. Or,     le crochet de Lie d'un \'el\'ement appartenant \`a la diff\'erence  $\mathcal G \setminus  heis$ avec un \'el\'ement de
    $heis \setminus \RR   X'$ n'est jamais contenu dans $\RR   X'$. Notre alg\`ebre ne peut donc pas \^etre nilpotente et de dimension $3$.

    Il vient que $\overline{\Gamma}$  est une copie de  $\RR^2 \subset G$ engendr\'ee par  $\RR X'$ et par n'importe quel  sous-groupe
    \`a un param\`etre $l$ contenu dans $\RR^2$ et transverse \`a $\RR   X'$. On conclut alors comme dans le cas pr\'ec\'edent. Plus pr\'ecis\'ement, le quotient de l'image de
    l'application d\'eveloppante par $\overline{\Gamma}$ s'identifie au quotient d'un ouvert de $SOL/I$ par le feuilletage trivial engendr\'e par $X'$. Or, ce quotient est s\'epar\'e et
    non compact et ne peut \^etre l'image continue de $M$ : absurde.

   Nous venons de d\'emontrer que (\`a indice fini  pr\`es)  $\Gamma \subset Heis$ et donc que (\`a rev\^etement fini  pr\`es) $M$ est un quotient de Heis
   par un r\'eseau cocompact. 
    \end{demonstration}
    
    \subsection{G non resoluble}
    
    Il reste \`a r\'egler le cas de la g\'eom\'etrie produit apparue au point  $(2)$ de la proposition~\ref{semi-simple algebrique}. Dans ce cas, $G= \RR \times  \widetilde{SL(2, \RR)}$ et $I \subset  \widetilde{SL(2, \RR)}$ est un groupe \`a un param\`etre semi-simple.
    
    \begin{proposition}   \label{non realisation compacte}
    Il n'existe pas de r\'ealisation compacte de $(G,G/I)$.
    \end{proposition}
    
    \begin{demonstration}
     Il s'agit d'un cas simple d'application de la m\'ethode d\'evelopp\'ee dans~\cite{Benoist-Labourie} (voir \'egalement la section 2.2.2 du  rapport de survol~\cite{Labourie}) dont nous
    pr\'esentons bri\`evement le principe. 
    
    La projection naturelle $G \to G/I$ munit  $G/I$ d'un fibr\'e en droites r\'eelles $G$-invariant, qui est donc bien d\'efini sur $M$. La forme lin\'eaire sur $\mathcal G$, qui s'annule sur le centre $\RR$ et qui co\"{\i}ncide sur $sl(2, \RR)$ avec la forme  duale (par rapport \`a la forme de Killing) d'un  g\'en\'erateur de $\mathcal I$,  \'equipe ce fibr\'e d'une connexion 
    $G$-invariante dont la forme de courbure $\omega$ est une forme volume en restriction aux orbites de $SL(2, \RR)$. La $2$-forme $\omega$ est alors bien d\'efinie sur $M$, o\`u elle s'annule sur $X$ et est non d\'eg\'en\'er\'ee en restriction aux feuilles de $X^{\bot}$.
    
    Par ailleurs, le fibr\'e en droites pr\'ec\'edent admettant une section jamais nulle (eventuellement sur un rev\^etement double de $M$), la forme de courbure $\omega$
    est exacte, \'egale \`a la diff\'erentielle de la $1$-forme de connexion $\omega_{1}$ associ\'ee \`a une section.
    
    Consid\'erons \`a pr\'esent la $1$-forme diff\'erentielle $\mu$ sur $M$ qui s'annule sur $X^{\bot}$ et telle que $\mu (X)=1$ (il s'agit de la forme duale
    de $X$ par rapport \`a $g$).  Comme le  flot de $X$ pr\'eserve le champ de plans int\'egrable $X^{\bot}$, la $1$-forme $\mu$ est ferm\'ee.
    
    La forme diff\'erentielle $\omega \wedge \mu$ est  une forme volume sur $M$. Par ailleurs, cette forme est exacte, \'egale \`a $d(\omega_{1} \wedge  \mu)$. Si $M$ est compacte, ceci
    est impossible.
  \end{demonstration}

\end{document}